\documentclass[hidelinks]{siamart251216}
\allowdisplaybreaks[2]

\usepackage{lipsum}
\usepackage{amsfonts}
\usepackage{graphicx}
\usepackage{epstopdf}
\usepackage{algorithmic}
\usepackage{hyperref}
\usepackage{enumitem}
\usepackage{multirow}
\usepackage{booktabs}
\usepackage{amssymb}
\usepackage{lmodern}
\usepackage[caption=false]{subfig}
\usepackage[T1]{fontenc}  
\setlist[itemize]{left=0pt}
\ifpdf
  \DeclareGraphicsExtensions{.eps,.pdf,.png,.jpg}
\else
  \DeclareGraphicsExtensions{.eps}
\fi

\newsiamthm{rek}{Remark}
\newsiamthm{alg}{Algorithm}
\newsiamremark{hypothesis}{Hypothesis}
\crefname{hypothesis}{Hypothesis}{Hypotheses}
\newsiamthm{claim}{Claim}
\newsiamremark{fact}{Fact}
\crefname{fact}{Fact}{Facts}

\headers{High-Order Area-Preserving Runge--Kutta Methods}{LEI CUI  AND GENGEN Zhang}

\title{High-Order Area-Preserving Runge--Kutta Methods for Surface Diffusion Flow of Curves \thanks{Submitted to the editors DATE.
\funding{This work was supported by the National Natural Science Foundation of China (No. 12261103), the Yunnan Fundamental Research Project (No. 202301AT070117),  and Postgraduate Scientific Research Innovation Fund Project of Yunnan University (No. KC-262615351).}}}

\author{Lei Cui\thanks{School of Mathematics and Statistics, Yunnan University, Kunming 650504,  China   (\email{isleicui@163.com}).}
\and Gengen Zhang\thanks{Corresponding author. School of Mathematics and Statistics, Yunnan University, Kunming 650504,  China (\email{zhanggen036@ynu.edu.cn}).}}

\usepackage{amsopn}

\ifpdf
\hypersetup{
  pdftitle={High-Order Area-Preserving Runge--Kutta Methods for Surface Diffusion Flow of Curves},
  pdfauthor={Lei Cui And Gengen Zhang}
}  

\fi   

\begin{document}
\maketitle

\begin{abstract} 
We propose a class of high-order area-preserving Runge--Kutta (RK) methods for  the surface diffusion flow of curves, a curvature-driven geometric flow that decreases the perimeter while conserving the enclosed area.  By  coupling the parametric finite element approximation of the classical BGN formulation with symplectic RK time discretisation, we rigorously prove that the fully discrete scheme unconditionally conserves the enclosed area of curves. This conservation arises solely from the intrinsic symplectic structure of the RK methods, without introducing any additional constraints, auxiliary variables, or Lagrange multipliers to enforce structure preservation. To attain high temporal accuracy, we adopt $q$-stage Gauss RK methods, which are symplectic and yield order $2q$ (e.g., $q=1,2,3$). The resulting fully implicit, strongly nonlinear system is solved efficiently via an equivalent discrete formulation combined with a monolithic damped Newton iteration. Extensive numerical tests on various evolving curves validate the theoretical findings, confirming exact area preservation, high-order temporal convergence, and favorable mesh quality, while the Newton iteration count rapidly decreases to fewer than 8 per step, demonstrating the robustness of the tailored solver.
\end{abstract}


\section{Introduction}

Surface diffusion is a mass‑transport mechanism governing solid‑material interface evolution via atomic migration along surfaces, causing gradual shape changes. Studied in both two‑ and three‑dimensional settings, it plays a central role in epitaxial growth \cite{1992PRL}, thermal grooving \cite{1957JAPMullins}, and solid‑state dewetting \cite{2017JCPBao,20SIAMJMAJW,15PRBWangY,20SISCZQ}. This work focuses exclusively on the two‑dimensional case.

In two dimensions, Mullins \cite{1957JAPMullins} modeled surface diffusion as the evolution of smooth closed curves $\Gamma(t)\subset\mathbb{R}^2$, $t\in[0,T]$, via
\begin{equation}\label{eq1.1}
	\mathcal{V}=\partial_{ss}\kappa\,\mathbf{n},
\end{equation}
where $\mathcal{V}$ is the velocity vector, $s$ the arclength parameter, $\kappa$ the curvature, and $\mathbf{n}$ the outward unit normal. This flow is the $H^{-1}$‑gradient flow of the perimeter functional \cite{2020HNABGN}. Let $L(t)$ and $A(t)$ denote perimeter and enclosed area, their temporal variations \cite{2005Acta} are
\begin{subequations}
	\begin{align}
		\frac{\mathrm{d}L(t)}{\mathrm{d}t}
		&=\int_{\Gamma(t)}\kappa(\mathcal{V}\cdot\mathbf{n})\,\mathrm{d}s
		=-\int_{\Gamma(t)}|\partial_s\kappa|^2\,\mathrm{d}s\leq 0,\quad t\geq 0,\\[1mm]
		\frac{\mathrm{d}A(t)}{\mathrm{d}t}
		&=\int_{\Gamma(t)}(\mathcal{V}\cdot\mathbf{n})\,\mathrm{d}s
		=\int_{\Gamma(t)}\partial_{ss}\kappa\,\mathrm{d}s=0,\quad t\geq 0.
	\end{align}
\end{subequations}
Thus,  the surface diffusion flow of closed curves decreases the perimeter while exactly preserving the enclosed area. This intrinsic area conservation is a geometric property that any robust numerical scheme should ideally mimic.

A pivotal reformulation of \eqref{eq1.1} was proposed by Barrett, Garcke, and N\"urnberg (BGN) \cite{2007JCPBGN}. Let $\mathbf{X}(\cdot,t):\mathbb{I}=[0,1]\rightarrow\mathbb{R}^2$ be a parametrization of the evolving closed curve $\Gamma(t)$. Under this reformulation, the governing equations become the coupled system
\begin{equation}\label{BGN}
	\begin{aligned}
		\partial_t\mathbf{X}\cdot\mathbf{n}&=\partial_{ss}\kappa,\\
		\kappa\mathbf{n}&=-\partial_{ss}\mathbf{X}.
	\end{aligned}
\end{equation}
The parametric finite element method (PFEM)  based on the BGN formulation exhibits favorable tangential mesh redistribution properties and energy stability \cite{2007JCPBGN,2020HNABGN}.  Bao and Zhao \cite{21SINUMBao} developed a first‑order structure‑preserving PFEM  for surface diffusion. Their approach uses a time‑weighted normal defined from two consecutive discrete interfaces to exactly conserve enclosed area or volume while ensuring energy dissipation. This idea has since been extended to several other geometric evolution equations \cite{23NMPDEBao,2024NMBao,2025MCBao}. However, these methods are limited to first‑order temporal accuracy.

The pursuit of higher temporal accuracy has motivated high‑order numerical methods for geometric flows. For spatial discretization, Garcke et al. \cite{2025JCPG} constructed isoparametric finite element methods of arbitrary polynomial degree for mean curvature flow and surface diffusion. For time integration, Jiang et al. proposed several higher‑order BGN‑based schemes, including a second‑order Crank–Nicolson–leapfrog method \cite{2024JCPCN},  higher‑order BDF schemes \cite{2024SISCBDFk}, and a predictor‑corrector method \cite{2025JCPPC}. In parallel, high‑order or optimal‑order error estimates have been established for various geometric flows and related problems \cite{2022NMLi,2019NMLubich,2021NMLubich,26LiarXiv}. Yet, these high‑order schemes do not generally preserve the area‑conservation property, highlighting a fundamental tension between high‑order accuracy and exact structure preservation.

This tension has led to growing interest in high‑order structure‑preserving discretizations. Garcke et al. \cite{2025SISCLagrange} introduced Lagrange multiplier formulations that incorporate area/volume preservation and energy dissipation into higher‑order time discretizations. Pan et al. \cite{2026Pan} constructed second‑order structure‑preserving PFEM  by exploiting exact geometric identities along a quadratic temporal path, while Zhang et al. \cite{2026ZGH} developed arbitrary‑order schemes based on a dual‑MDR formulation combined with continuous Petrov–Galerkin time discretization. Further related contributions can be found in \cite{2026arXivDuan,2025JCPGaoLi,2025JCPGaoY,2025CiCPLiM}. Nevertheless, constructing a general framework that simultaneously delivers \emph{arbitrarily high‑order temporal accuracy} and \emph{unconditional area preservation} for the surface diffusion flow of curves remains an important open challenge.

In contrast, Runge–Kutta (RK) time discretization offers distinct advantages in simultaneously achieving high‑order accuracy and structure preservation. On the one hand, RK methods are inherently self‑starting, eliminating the need for historical data storage. This feature yields several concrete benefits: (i) natural suitability for adaptive step‑size control via embedded error estimators; (ii) a straightforward path to high‑order extension by adding stages, avoiding the zero‑stability constraints typical of linear multistep methods; and (iii) favorable stability properties, such as A‑stability and L‑stability \cite{Ascher,1996Hairer}. On the other hand, RK methods also play a key role in structure preservation. For instance, symplectic RK methods preserve the symplectic structure of Hamiltonian systems \cite{2010FengKang,Hairer} and have been widely used to maintain conserved quantities in various conservative systems \cite{2023LiXJCP,2023AdvZhang}. Furthermore, algebraically stable RK methods are commonly employed to preserve dissipativity in dissipative systems \cite{2019SISCAkrivis,26BaoarXiv,2020JCPGYZ}.

We exploit this dual capability, high‑order accuracy and intrinsic structure preservation, in this work. We propose a class of high‑order area‑preserving RK  methods for the parametric finite element approximation of surface diffusion flow of curves. Starting from the semi‑discretization of the classical BGN formulation, we recast it into a differential‑algebraic equations framework, which facilitates a transparent exposition of the temporal discretization. After discretizing in time with symplectic RK methods, we rigorously prove that the resulting scheme \emph{unconditionally conserves the enclosed area} of curves. Crucially, this conservation property stems entirely from the intrinsic symplectic structure of the RK methods, without requiring any additional constraints. To achieve the desired high‑order temporal accuracy, we select Gauss RK methods that satisfy the symplectic condition, thereby attaining an order of $2q$ for $q$ stages. For the fully implicit, strongly nonlinear system with coupled stages, we design a dedicated solver based on a monolithic damped Newton iteration.  The main contributions of this paper are summarized as follows:
\begin{enumerate}[label=\(\bullet\)]
	\item \textbf{Unconditional area preservation.} 
	We prove that symplectic RK schemes unconditionally conserve the enclosed area, a property dictated solely by the intrinsic symplectic structure, requiring no additional constraints,  auxiliary variables or  Lagrange multipliers to impose structure preservation.

	\item \textbf{Arbitrarily high‑order temporal accuracy.} 
	Adopting $q$-stage Gauss RK methods (symplectic, order $2q$) yields a family of schemes with arbitrarily high convergence rates (second, fourth, and sixth), providing a novel and efficient pathway to high‑order temporal algorithms within the BGN formulation while retaining exact area preservation.
	
	\item \textbf{Efficient solution of the fully implicit system.} 
A tailored monolithic damped Newton iteration, built upon an equivalent discrete formulation, efficiently solves the coupled nonlinear system; numerical experiments demonstrate a rapid decrease to under 8 Newton iterations per step, confirming the robustness and efficiency of the tailored solver.
\end{enumerate}

The remainder of this paper is organized as follows. Section \ref{Section2} applies PFEM to spatially discretize the classical BGN formulation \eqref{BGN} and presents the corresponding DAE form. Section \ref{Section3} introduces symplectic RK methods for temporal discretization and rigorously proves that the fully discrete scheme unconditionally preserves the area of closed curves. Section \ref{Section4} constructs an equivalent formulation for the fully implicit RK discretization and develops a monolithic damped Newton iteration scheme to solve the coupled nonlinear system efficiently. Numerical experiments are presented in Section \ref{Section5}, where we verify the area-preserving property and high-order convergence rates, and demonstrate the rapid drop in Newton iterations per step, thereby showcasing the tailored solver's efficiency. Finally, Section \ref{Section6} summarizes the findings of this work.

\section{Weak formulation and spatial discretization}\label{Section2}
This section employs a PFEM for the spatial discretization of the classical BGN formulation \eqref{BGN} of surface diffusion flow.  
Let $\mathbb{I}=[0,1]$ be partitioned into $N\ge 3$ subintervals $I_j = [\rho_{j-1},\rho_j]$, $j=1,\dots,N$, with mesh size $h=\max_{1\leq j\leq N}|\rho_j-\rho_{j-1}|$.  
We define the function space on $\Gamma(t)=\mathbf{X}(\mathbb{I})$ as
\begin{equation}
	L^2(\mathbb{I}) := \Bigl\{ u : \mathbb{I} \to \mathbb{R} \;\Big|\;  \int_{\Gamma(t)} \left|u(s) \right|^2\,\mathrm{d}s=\int_{\mathbb{I}}^{}\left|u\circ\mathbf{X}(\rho) \right|^2\left|\partial_\rho\mathbf{X} \right|\,\mathrm{d}\rho<+\infty   \Bigr\},
\end{equation}
equipped with the inner product
\begin{equation}\label{eq2.2}
	\left(u,v\right)_{\Gamma}:=\int_{\Gamma(t)}^{}u(s)v(s)\,\mathrm{d}s=\int_{\mathbb{I}}^{}\left(u\circ\mathbf{X}(\rho) \right) \left(v\circ\mathbf{X}(\rho) \right)\left|\partial_\rho\mathbf{X} \right|\,\mathrm{d}\rho, \ \forall\,u,v\in L^2(\mathbb{I}).
\end{equation}
The Sobolev space on the reference domain is defined as
\begin{equation}
	H^1(\mathbb{I}):=\bigl\{u:\mathbb{I}\rightarrow\mathbb{R}\;\big|\;u\in L^2(\mathbb{I}),\;\partial_\rho u\in L^2(\mathbb{I})\bigr\}.
\end{equation}
We introduce the periodic piecewise‑linear finite element space
\begin{equation}\nonumber
	V^h := \bigl\{ u \in C(\mathbb{I}) : u|_{I_j} \text{ is linear},\; \forall j=1,\dots,N,\; u(\rho_0)=u(\rho_N) \bigr\} \subseteq   H^1(\mathbb{I}).
\end{equation}
For a polygonal curve $\Gamma^h = \mathbf{X}^h(\mathbb{I}) \in [V^h]^2$, the inner product \eqref{eq2.2} is replaced by its mass‑lumped counterpart obtained via the composite trapezoidal rule:
\begin{equation}\label{lumped_inner}
	(u,v)_{\Gamma^h}^h := \frac12 \sum_{j=1}^N |\mathbf{h}_j| \bigl[ (u\cdot v)(\rho_j^-) + (u\cdot v)(\rho_{j-1}^+) \bigr],
\end{equation}
where $\mathbf{h}_j := \mathbf{X}^h(\rho_j,t) - \mathbf{X}^h(\rho_{j-1},t)$,  and $u(\rho_j^{\pm}) = \lim_{\rho\to\rho_j^{\pm}} u(\rho)$ denotes the one‑sided limits for functions that may be discontinuous at the nodes.

Using integration by parts and periodicity, the weak form of \eqref{BGN} reads: given an initial closed curve $\Gamma^h(0)=\mathbf{X}^h(\cdot,0)\in[V^h]^2$ parametrized clockwise,  find the evolving curve $\mathbf{X}^h(\cdot,t)\in[V^h]^2$ and its curvature $\kappa^h(\cdot,t)\in V^h$ such that
\begin{equation}\label{eq2.5}
	\begin{aligned}
		\bigl(\partial_t\mathbf{X}^h \cdot \mathbf{n}^h, \; \varphi^h \bigr)_{\Gamma^h}^h + \bigl( \partial_s\kappa^h,\; \partial_s\varphi^h \bigr)_{\Gamma^h} &= 0, && \forall\,\varphi^h \in V^h, \\[2mm]
		\bigl( \kappa^h,\; \mathbf{n}^h \cdot \boldsymbol{\omega}^h \bigr)_{\Gamma^h}^h - \bigl( \partial_s\mathbf{X}^h,\; \partial_s\boldsymbol{\omega}^h \bigr)_{\Gamma^h} &= 0, && \forall\,\boldsymbol{\omega}^h \in [V^h]^2 .
	\end{aligned}
\end{equation}
Here the outward unit normal $\mathbf{n}^h$ is piecewise constant on each edge,
\begin{equation}\nonumber
	\mathbf{n}^h|_{I_j} = - \frac{\mathbf{h}_j^\perp}{|\mathbf{h}_j|}, \quad j=1,\cdots,N,
\end{equation}
where $(\cdot)^\perp$ denotes clockwise rotation by $\pi/2$.   The arc‑length derivative $\partial_s$ on each segment $I_j$ is approximated by
\begin{equation}\nonumber
	\partial_s f|_{I_j}=\frac{\partial_\rho f}{|\partial_\rho\mathbf{X}^h|}|_{I_j}=\frac{(\rho_j-\rho_{j-1})\partial_\rho f|_{I_j}}{|\mathbf{h}_j|}.
\end{equation}
As established in \cite{2007JCPBGN,2020HNABGN}, the semi‑discrete scheme \eqref{eq2.5} yields asymptotically equidistributed meshes, provided the vertices are not locally parallel.

To recast \eqref{eq2.5} as a differential‑algebraic equations, let $\{\varphi_i\}_{i=1}^N$ be the standard periodic nodal basis of $V^h$, satisfying $\varphi_i(\rho_j)=\delta_{ij}$ (the Kronecker delta).  
Expanding the discrete solutions as
\begin{equation}\nonumber
	\mathbf X^h(\rho,t)=\sum_{i=1}^N\mathbf X_i(t)\varphi_i(\rho),
	\qquad
	\kappa^h(\rho,t)=\sum_{i=1}^N\kappa_i(t)\varphi_i(\rho),
\end{equation}
with $\mathbf X_i(t)=(x_i(t),y_i(t))^T$, we collect the coefficients into vectors
\begin{equation}\nonumber
	\mathbf U=
	\begin{pmatrix}
		\mathbf x_1\\
		\mathbf x_2
	\end{pmatrix}\in\mathbb R^{2N},
	\qquad
	\mathbf K=(\kappa_1,\dots,\kappa_N)^T\in\mathbb R^N,
\end{equation}
where $\mathbf x_1=(x_1,\dots,x_N)^T$ and $\mathbf x_2=(y_1,\dots,y_N)^T$.

We define the stiffness matrix $P(\mathbf U)\in\mathbb R^{N\times N}$ and the normal coupling matrices $M_\ell(\mathbf U)\in\mathbb R^{N\times N}$ ($\ell=1,2$) by
\begin{equation}\nonumber
	(P(\mathbf U))_{ij}
	=(\partial_s\varphi_i,\partial_s\varphi_j)_{\Gamma^h},\quad
	(M_\ell(\mathbf U))_{ij}
	=(n_\ell^h\varphi_i,\varphi_j)^h_{\Gamma^h},
	\quad 1\le i,j\le N,
\end{equation}
where $n_\ell^h$ is the $\ell$-th component of $\mathbf n^h$.  
Owing to the mass‑lumping \eqref{lumped_inner}, each $M_\ell(\mathbf U)$ is diagonal with entries
\begin{equation}\nonumber
	(M_\ell(\mathbf U))_{jj}
	=\frac12\Bigl(
	|\mathbf h_j|\,n_{j,\ell}^h
	+|\mathbf h_{j+1}|\,n_{j+1,\ell}^h
	\Bigr),
\end{equation}
where $n_{j,\ell}^h$ denotes the $\ell$-th component of the outward unit normal  $\mathbf n^h$ on edge $I_j$, and indices are taken periodically, in particular, $\mathbf h_{N+1}=\mathbf h_1$ and \(n_{N+1,\ell}^h=n_{1,\ell}^h\).

Choosing test functions $\varphi^h=\varphi_i$ and $\boldsymbol{\omega}^h=\varphi_i\mathbf{e}_\ell$ ($\ell=1,2$) in \eqref{eq2.5} yields the algebraic system
\begin{equation}\label{eq2.6}
	\begin{cases}
		M_1(\mathbf U)\dot{\mathbf x}_1
		+M_2(\mathbf U)\dot{\mathbf x}_2
		+P(\mathbf U)\mathbf K=0,\\
		M_1(\mathbf U)\mathbf K-P(\mathbf U)\mathbf x_1=0,\\
		M_2(\mathbf U)\mathbf K-P(\mathbf U)\mathbf x_2=0,
	\end{cases}
\end{equation}
where $\dot{z}=\partial_t z$, and \(\mathbf{e}_\ell\) denotes the \(\ell\)-th Cartesian unit vector in \(\mathbb R^2\).  
Introducing the block matrices
\[
\mathbb M(\mathbf U)=
\begin{pmatrix}
	M_1(\mathbf U)&M_2(\mathbf U)
\end{pmatrix}\in\mathbb R^{N\times2N},
\quad
\mathbb P(\mathbf U)=
\begin{pmatrix}
	P(\mathbf U)&0\\
	0&P(\mathbf U)
\end{pmatrix}\in\mathbb R^{2N\times2N},
\]
we rewrite \eqref{eq2.6} compactly as
\begin{equation}\label{eq2.7}
	\begin{cases}
		\mathbb M(\mathbf U)\dot{\mathbf U}
		+P(\mathbf U)\mathbf K=0,\\
		\mathbb M(\mathbf U)^T\mathbf K
		-\mathbb P(\mathbf U)\mathbf U=0.
	\end{cases}
\end{equation}
This DAE system serves as the foundation for the time‑stepping schemes discussed in the next section.

\section{Temporal discretization}\label{Section3}

RK methods are widely used for the numerical integration of differential equations due to their favorable stability properties and high‑order accuracy \cite{Ascher,1996Hairer}. In this section, we apply symplectic RK methods to discretize in time the semi‑discrete system \eqref{eq2.7}. A distinguishing feature of the resulting scheme is its unconditional area‑conservation property for the surface diffusion flow of closed curves.

\subsection{Symplectic Runge--Kutta methods}
Let $\tau$ denote the time step size  and set $t_m=m\tau$. We introduce the intermediate-stage times as $t_{m,i}=t_m+c_i\tau$  and approximate the solution triples $(\mathbf{U},\partial_t\mathbf{U},\mathbf{K})$ at these intermediate stages by $(\mathbf{U}^{m,i}, \mathbf{V}^{m,i}, \mathbf{K}^{m,i})$. Applying a $q$-stage symplectic RK method to \eqref{eq2.7} yields the stage equations
\begin{equation}\label{eq3.1}
	\mathbf{U}^{m,i}=\mathbf{U}^{m}+\tau\sum_{j=1}^{q}a_{ij}\mathbf{V}^{m,j},\qquad i=1,\dots,q,
\end{equation}
where the stage vectors satisfy the algebraic  constraints
\begin{equation}
	\begin{cases}
		\mathbb{M}(\mathbf{U}^{m,i})\mathbf{V}^{m,i}
		+\mathbb{P}(\mathbf{U}^{m,i})\mathbf{K}^{m,i}=0,\\[1mm]
		\mathbb{M}(\mathbf{U}^{m,i})^{\!T}\mathbf{K}^{m,i}
		-\mathbb{P}(\mathbf{U}^{m,i})\mathbf{U}^{m,i}=0.
	\end{cases}
\end{equation}
The numerical solution is advanced by
\begin{equation}\label{eq3.2}
	\mathbf{U}^{m+1}=\mathbf{U}^{m}+\tau\sum_{i=1}^{q}b_{i}\mathbf{V}^{m,i},
\end{equation}
and the curvature vector $\mathbf{K}^{m+1}$ is determined by the discrete constraint
\begin{equation}\label{eqK}
	\mathbb{M}(\mathbf{U}^{m+1})^{\!T}\mathbf{K}^{m+1}=\mathbb{P}(\mathbf{U}^{m+1})\mathbf{U}^{m+1}.
\end{equation}
When the coefficients satisfy $b_i=a_{qi}$ for $i=1,\dots,q$, the update simplifies to $\mathbf{U}^{m+1}=\mathbf{U}^{m,q}$ and $\mathbf{K}^{m+1}=\mathbf{K}^{m,q}$.

The coefficients $(\mathbf A,\mathbf{b},\mathbf{c})$ of the symplectic RK scheme are displayed in the Butcher tableau
\begin{equation}\label{eq3.3}
	\begin{array}{c|c}
		\mathbf{c} & \mathbf A \\ \hline
		& \mathbf{b}^{T}
	\end{array}
	=
	\begin{array}{c|cccc}
		c_1 & a_{11} & a_{12} & \cdots & a_{1q} \\
		c_2 & a_{21} & a_{22} & \cdots & a_{2q} \\
		\vdots & \vdots & \vdots & \ddots & \vdots \\
		c_q & a_{q1} & a_{q2} & \cdots & a_{qq} \\ \hline
		& b_1 & b_2 & \cdots & b_q
	\end{array},
\end{equation}
which obeys the consistency condition $c_i=\sum_{j=1}^{q}a_{ij}$ and the symplecticity condition \cite{2010FengKang,Hairer}
\begin{equation}\label{Symplectic}
	b_i a_{ij} + b_j a_{ji} = b_i b_j , \qquad i,j=1,\dots,q.
\end{equation}

To analyze area conservation, we rewrite the above  coupled nonlinear algebraic system into an equivalent fully discrete variational form. Let $\mathbf{X}^{m,i}$, $\boldsymbol{\Phi}^{m,i}$, and $\kappa^{m,i}$ denote the approximations of $\mathbf{X}(\cdot,t_{m,i})$, $\partial_t\mathbf{X}(\cdot,t_{m,i})$, and $\kappa(\cdot,t_{m,i})$, respectively. The stage positions are given by
\begin{equation}\label{eq3.6}
	\mathbf{X}^{m,i}=\mathbf{X}^{m}+\tau\sum_{j=1}^{q}a_{ij}\boldsymbol{\Phi}^{m,j},\qquad i=1,\dots,q.
\end{equation}
We then seek $\mathbf{X}^{m,i},\boldsymbol{\Phi}^{m,i}\in[V^h]^2$ and $\kappa^{m,i}\in V^h$ such that \eqref{eq3.6} holds together with the weak formulation: for each $i=1,\dots,q$,
\begin{equation}\label{eq3.7}
	\begin{aligned}
		\bigl( \boldsymbol{\Phi}^{m,i} \cdot \mathbf{n}^{m,i}, \, \varphi^h \bigr)_{\Gamma^{m,i}}^h
		+ \bigl( \partial_s^{m,i}\kappa^{m,i},\, \partial_s^{m,i}\varphi^h \bigr)_{\Gamma^{m,i}} &=0,
		&&\forall\, \varphi^h \in V^h, \\[1mm]
		\bigl( \kappa^{m,i},\,\mathbf{n}^{m,i} \cdot \boldsymbol{\omega}^h \bigr)_{\Gamma^{m,i}}^h
		- \bigl( \partial_s^{m,i}\mathbf{X}^{m,i},\, \partial_s^{m,i}\boldsymbol{\omega}^h \bigr)_{\Gamma^{m,i}} &=0,
		&&\forall\, \boldsymbol{\omega}^h \in [V^h]^2,
	\end{aligned}
\end{equation}
where $\Gamma^{m,i}:=\mathbf{X}^{m,i}(\mathbb{I})$, $\mathbf{n}^{m,i}$ is the outward unit normal on $\mathbf{X}^{m,i}$, and $\partial_s^{m,i}f=|\partial_\rho \mathbf{X}^{m,i}|^{-1}\partial_\rho f$. The update for the next time level is
\begin{equation}\label{eq3.8}
	\mathbf{X}^{m+1}=\mathbf{X}^{m}+\tau\sum_{i=1}^{q}b_{i}\boldsymbol{\Phi}^{m,i}.
\end{equation}
The curvature $\kappa^{m+1}$ satisfies the algebraic equation
\begin{equation}\label{eqka}
	\bigl( \kappa^{m+1},\mathbf{n}^{m+1} \cdot \boldsymbol{\omega}^h \bigr)_{\Gamma^{m+1}}^h
	- \bigl( \partial_s\mathbf{X}^{m+1}, \partial_s^{m+1}\boldsymbol{\omega}^h \bigr)_{\Gamma^{m+1}} =0,
	\qquad\forall\, \boldsymbol{\omega}^h \in [V^h]^2,
\end{equation}
which can be solved efficiently via a least‑squares procedure. 

\begin{rek}
	The matrix‑oriented scheme \eqref{eq3.1}--\eqref{eqK} and the variational scheme \eqref{eq3.6}--\eqref{eqka} are fully equivalent. The area‑preserving property of the latter is proved in the next subsection. Because the method is fully implicit, all intermediate stages are coupled, leading to a strongly nonlinear system. In Section \ref{Section4} we propose a tailored iterative solver to mitigate this cost.
\end{rek}

\subsection{Unconditional area‑preserving property}
We now establish that the symplectic RK scheme \eqref{eq3.6}--\eqref{eqka} exactly conserves the enclosed area. Let $A^{m}$ denote the discrete area of $\Gamma^m$:
\begin{equation}\label{eqAm}
	A^m:=\frac{1}{2}\int_{\Gamma^m}^{}\mathbf{X}^m\cdot\mathbf{n}^m\,\mathrm{d}s=-\frac{1}{2}\int_{\mathbb{I}}^{}\mathbf{X}^m\cdot(\partial_\rho \mathbf{X}^{m})^\perp \mathrm{d}\rho.
\end{equation}

\begin{theorem}[Unconditional area conservation]\label{thmAP}
	Let $(\mathbf{X}^{m+1},\kappa^{m+1})$ be the solution produced by the symplectic RK scheme \eqref{eq3.6}--\eqref{eqka}. Then
	\begin{equation}
		A^{m+1}=A^m=\cdots=A^0,\qquad m\ge 0.
	\end{equation}
\end{theorem}

\begin{proof}
	Substituting the update \eqref{eq3.8} into \eqref{eqAm} gives
	\begin{equation}\nonumber
		\begin{aligned}
			A^{m+1}
			&=-\frac{1}{2}\int_{\mathbb{I}}\Bigl(\mathbf{X}^{m}+\tau\sum_{i=1}^{q}b_{i}\boldsymbol{\Phi}^{m,i} \Bigr)
			\cdot\Bigl( \partial_\rho \mathbf{X}^{m}+\tau\sum_{i=1}^{q}b_{i}\partial_\rho\boldsymbol{\Phi}^{m,i} \Bigr)^\perp \,\mathrm{d}\rho \\[1mm]
			&=A^m-\frac{\tau}{2}\sum_{i=1}^{q}b_i\int_{\mathbb{I}}\mathbf{X}^{m}\cdot\bigl(\partial_\rho\boldsymbol{\Phi}^{m,i}\bigr)^\perp \,\mathrm{d}\rho
			-\frac{\tau}{2}\sum_{i=1}^{q}b_i\int_{\mathbb{I}}\boldsymbol{\Phi}^{m,i}\cdot\bigl(\partial_\rho\mathbf{X}^{m}\bigr)^\perp \,\mathrm{d}\rho \\[1mm]
			&\quad-\frac{\tau^2}{2}\sum_{i,j=1}^{q}b_ib_j\int_{\mathbb{I}}\boldsymbol{\Phi}^{m,i}\cdot\bigl(\partial_\rho\boldsymbol{\Phi}^{m,j}\bigr)^\perp \,\mathrm{d}\rho.
		\end{aligned}
	\end{equation}
	Using the identity $\mathbf{v}\cdot\mathbf{w}^\perp=-\mathbf{v}^\perp\cdot\mathbf{w}$ and periodicity, integration by parts yields
	\begin{equation}\label{eq3.13}
		\int_{\mathbb{I}}\mathbf{X}^{m}\cdot\bigl(\partial_\rho\boldsymbol{\Phi}^{m,i}\bigr)^\perp \,\mathrm{d}\rho
		=\int_{\mathbb{I}}\boldsymbol{\Phi}^{m,i}\cdot\bigl(\partial_\rho\mathbf{X}^{m}\bigr)^\perp \,\mathrm{d}\rho,
		\qquad i=1,\dots,q.
	\end{equation}
	Hence,
	\begin{equation}\nonumber
		A^{m+1}-A^{m}
		=-\tau\sum_{i=1}^{q}b_i\int_{\mathbb{I}}\boldsymbol{\Phi}^{m,i}\cdot\bigl(\partial_\rho\mathbf{X}^{m}\bigr)^\perp \,\mathrm{d}\rho
		-\frac{\tau^2}{2}\sum_{i,j=1}^{q}b_ib_j\int_{\mathbb{I}}\boldsymbol{\Phi}^{m,i}\cdot\bigl(\partial_\rho\boldsymbol{\Phi}^{m,j}\bigr)^\perp \,\mathrm{d}\rho.
	\end{equation}
Next, using the stage relation \eqref{eq3.6}, we replace $\mathbf{X}^{m}$ by $\mathbf{X}^{m,i}-\tau\sum\limits_{j=1}^{q}a_{ij}\boldsymbol{\Phi}^{m,j}$:
	\begin{equation}\label{eq3.14}
		\begin{aligned}
			A^{m+1}-A^{m}
			&=-\tau\sum_{i=1}^{q}b_i\int_{\mathbb{I}}\boldsymbol{\Phi}^{m,i}\cdot\bigl(\partial_\rho\mathbf{X}^{m,i}\bigr)^\perp \,\mathrm{d}\rho \\
			&\quad+\tau^2\sum_{i,j=1}^{q}b_ia_{ij}\int_{\mathbb{I}}\boldsymbol{\Phi}^{m,i}\cdot\bigl(\partial_\rho\boldsymbol{\Phi}^{m,j}\bigr)^\perp \,\mathrm{d}\rho \\
			&\quad-\frac{\tau^2}{2}\sum_{i,j=1}^{q}b_ib_j\int_{\mathbb{I}}\boldsymbol{\Phi}^{m,i}\cdot\bigl(\partial_\rho\boldsymbol{\Phi}^{m,j}\bigr)^\perp \,\mathrm{d}\rho.
		\end{aligned}
	\end{equation}

By the same integration‑by‑parts argument as in \eqref{eq3.13}, one finds the relation
$\int_{\mathbb{I}} \mathbf{\Phi}^{m,j}\cdot(\partial_\rho\mathbf{\Phi}^{m,i})^\perp\,\mathrm{d}\rho
= \int_{\mathbb{I}} \mathbf{\Phi}^{m,i}\cdot(\partial_\rho\mathbf{\Phi}^{m,j})^\perp\,\mathrm{d}\rho, $
which shows that 
\begin{equation}\label{eq3.15}
	\sum_{i,j=1}^{q}b_ia_{ij}\int_{\mathbb{I}}^{}\mathbf{\Phi}^{m,i}\cdot\left( \partial_\rho\mathbf{\Phi}^{m,j} \right)^\perp\mathrm{d}\rho=\sum_{i,j=1}^{q}b_ja_{ji}\int_{\mathbb{I}}^{}\mathbf{\Phi}^{m,i}\cdot\left( \partial_\rho\mathbf{\Phi}^{m,j} \right)^\perp \mathrm{d}\rho.
\end{equation}
Substituting \eqref{eq3.15} into \eqref{eq3.14} and invoking the symplectic condition \eqref{Symplectic} leads to
	\begin{equation}
		\begin{aligned}
			A^{m+1}-A^{m}
			&=-\tau\sum_{i=1}^{q}b_i\int_{\mathbb{I}}\boldsymbol{\Phi}^{m,i}\cdot\bigl(\partial_\rho\mathbf{X}^{m,i}\bigr)^\perp \,\mathrm{d}\rho \\
			&\quad+\frac{\tau^2}{2}\sum_{i,j=1}^{q}(b_ia_{ij}+b_ja_{ji}-b_ib_j)\int_{\mathbb{I}}\boldsymbol{\Phi}^{m,i}\cdot\bigl(\partial_\rho\boldsymbol{\Phi}^{m,j}\bigr)^\perp \,\mathrm{d}\rho \\[1mm]
			&=-\tau\sum_{i=1}^{q}b_i\int_{\mathbb{I}}\boldsymbol{\Phi}^{m,i}\cdot\bigl(\partial_\rho\mathbf{X}^{m,i}\bigr)^\perp \,\mathrm{d}\rho.
		\end{aligned}
	\end{equation}

	Finally, choosing the test function $\varphi^h\equiv 1$ in the first equation of \eqref{eq3.7} annihilates the second term (since $\partial_s 1=0$) and yields
\begin{equation}\nonumber
	-\int_{\mathbb{I}}^{}\mathbf{\Phi}^{m,i}\cdot\left( \partial_\rho\mathbf{X}^{m,i}\right)^\perp \mathrm{d}\rho=\int_{\Gamma^{m,i}}^{}\mathbf{\Phi}^{m,i}\cdot\mathbf{n}^{m,i} \,\mathrm{d}s=0,  \quad i=1,\dots, q,
\end{equation}
	which implies $A^{m+1}=A^m$. This completes the proof.
\end{proof}

In our numerical experiments, we primarily use Gauss RK methods, a family of high‑order symplectic schemes whose Butcher tableaux are listed in Appendix~\ref{appendixA}. A $q$-stage Gauss method achieves order $2q$ and is A‑stable \cite{Ascher,1996Hairer}, making it robust for stiff geometric evolution problems.

\begin{rek}
	Generic RK methods do not generally preserve geometric structures. Theorem \ref{thmAP} shows that RK schemes satisfying the symplectic condition \eqref{Symplectic} exactly conserve area for closed curves, a property crucial for accurately capturing the long‑time dynamics of complex shapes.
\end{rek}

\begin{rek}
	Unlike many classical structure‑preserving discretizations \cite{2025SISCLagrange,2021JCPJiangLi,2025CiCPLiM}, the present methods attain high-order accuracy without introducing  auxiliary variables or  Lagrange multipliers to impose structure preservation. Instead, they rely solely on the inherent symplecticity of the RK time integrator to achieve conservation.
	
\end{rek}

\section{Reformulation and iterative solution}\label{Section4}

Because the fully implicit RK discretization involves multiple stage velocities, outward unit normals, and curvatures, all dependent on the unknown stage curves,  the resulting fully discrete system is highly nonlinear. To solve it efficiently, we adopt the equivalent reformulation introduced by Barrett, Garcke, and N\"urnberg for fully implicit schemes for planar curves~\cite{2011BGN}, and combine it with a monolithic damped Newton method that handles all RK stages simultaneously.

\subsection{Reformulation of the first‑order fully implicit BGN scheme}

We first recall the first‑order fully implicit scheme associated with the semi‑discrete system~\eqref{eq2.5}: for $m\ge 0$, find $\bigl(\mathbf{X}^{m+1},\kappa^{m+1}\bigr)\in[V^h]^2\times V^h$ such that
\begin{equation}\label{eq4.1}
	\begin{aligned}
		\Bigl( \frac{\mathbf{X}^{m+1}-\mathbf{X}^m}{\tau},\; \varphi^h \mathbf{n}^{m+1} \Bigr)_{\Gamma^{m+1}}^h + \bigl( \partial_s^{m+1}\kappa^{m+1},\; \partial_s^{m+1}\varphi^h \bigr)_{\Gamma^{m+1}}&=0, \\[2mm]
		\bigl( \kappa^{m+1},\; \mathbf{n}^{m+1} \cdot \boldsymbol{\omega}^h \bigr)_{\Gamma^{m+1}}^h - \bigl( \partial_s^{m+1}\mathbf{X}^{m+1},\; \partial_s^{m+1}\boldsymbol{\omega}^h \bigr)_{\Gamma^{m+1}}&=0,
	\end{aligned}
\end{equation}
for all $\varphi^h \in V^h$ and $\boldsymbol{\omega}^h \in [V^h]^2$, where $\mathbf{n}^{m+1}$ is the unit outward normal on $\mathbf{X}^{m+1}$ and $\partial_s^{m+1}f=|\partial_\rho \mathbf{X}^{m+1}|^{-1}\partial_\rho f$.

Applying Newton's method directly to~\eqref{eq4.1} typically yields a singular Jacobian matrix and fails to converge. Fixed‑point iteration, on the other hand, suffers from prohibitively slow convergence. To overcome these difficulties, Barrett, Garcke, and N\"urnberg~\cite{2011BGN} proposed an equivalent formulation: find $\bigl(\mathbf{X}^{m+1},\kappa^{m+1}\bigr)\in[V^h]^2\times V^h$ such that
\begin{equation}\label{eq4.2}
	\begin{aligned}
		\Bigl( \frac{\mathbf{X}^{m+1}-\mathbf{X}^m}{\tau},\;  -(\partial_\rho\mathbf{X}^{m+1})^\perp\varphi^h \Bigr)_{\mathbb{I}}^h + \frac{1}{L^{m+1}}\bigl( \partial_\rho\kappa^{m+1},\; \partial_\rho\varphi^h \bigr)_{\mathbb{I}}&=0, \\[2mm]
		\bigl( \kappa^{m+1},\; -(\partial_\rho\mathbf{X}^{m+1})^\perp \cdot \boldsymbol{\omega}^h \bigr)_{\mathbb{I}}^h - \frac{1}{L^{m+1}}\bigl( \partial_\rho\mathbf{X}^{m+1},\; \partial_\rho\boldsymbol{\omega}^h \bigr)_{\mathbb{I}}&=0,
	\end{aligned}
\end{equation}
for all $\varphi^h \in V^h$ and $\boldsymbol{\omega}^h \in [V^h]^2$. Here $L^{m+1}=\int_{\mathbb{I}}|\partial_\rho\mathbf{X}^{m+1}|\,\mathrm{d}\rho$ denotes the perimeter, and the mass‑lumped inner product on $\mathbb{I}$ is defined by
\begin{equation}\label{Newproduct}
	(u,v)^h_{\mathbb I}
	=
	\frac12\sum_{j=1}^N\Delta\rho_j
	\left[
	(u\cdot v)(\rho_j^-)
	+
	(u\cdot v)(\rho_{j-1}^+)
	\right],
\end{equation}
for all $u,v \in L^2(\mathbb{I})$.

The equivalence between the reformulation \eqref{eq4.2} and the original nonlinear system~\eqref{eq4.1} relies on mild assumptions. Moreover, any solution of~\eqref{eq4.2} satisfies a strict mesh‑equidistribution property and automatically fulfills the original scheme~\eqref{eq4.1}. Thus the reformulation can be viewed as a strictly equidistributed solution branch of the original problem; details are available in~\cite{2011BGN}.

\subsection{Reformulation of the fully implicit RK discretization}

We now extend the above reformulation to the intermediate‑stage system~\eqref{eq3.7} of the fully implicit RK discretization and rigorously establish its relation to the original stage equations.

Let $\mathbf A=(a_{ij})_{i,j=1}^q$ be the coefficient matrix of the Butcher tableau, assumed nonsingular, and write $\mathbf A^{-1}=(\widehat{a}_{ij})_{i,j=1}^q$. From the stage equations~\eqref{eq3.6}, the stage velocities admit the representation
\begin{equation}
	\mathbf{\Phi}^{m,i}=\frac{1}{\tau}\sum_{j=1}^{q}\widehat{a}_{ij}(\mathbf{X}^{m,j}-\mathbf{X}^m),\qquad i=1,\dots,q.
\end{equation}
For each stage, let $L^{m,i}=\int_{\mathbb{I}}|\partial_\rho\mathbf{X}^{m,i}|\,\mathrm{d}\rho$ denote the perimeter of the curve $\Gamma^{m,i}$. Applying the reformulation~\eqref{eq4.2} to every RK stage yields the coupled nonlinear system: for each $i=1,\dots,q$, find $\mathbf X^{m,i},\mathbf\Phi^{m,i}\in[V^h]^2$ and $\kappa^{m,i}\in V^h$ satisfying
\begin{equation}\label{eq4.4}
	\begin{aligned}
		\Bigl( \mathbf{\Phi}^{m,i},\;  -(\partial_\rho\mathbf{X}^{m,i})^\perp\varphi^h \Bigr)_{\mathbb{I}}^h + \frac{1}{L^{m,i}}\bigl( \partial_\rho\kappa^{m,i},\; \partial_\rho\varphi^h \bigr)_{\mathbb{I}}&=0,\quad \forall\,\varphi^h \in V^h, \\[2mm]
		\bigl( \kappa^{m,i},\; -(\partial_\rho\mathbf{X}^{m,i})^\perp \cdot \boldsymbol{\omega}^h \bigr)_{\mathbb{I}}^h - \frac{1}{L^{m,i}}\bigl( \partial_\rho\mathbf{X}^{m,i},\; \partial_\rho\boldsymbol{\omega}^h \bigr)_{\mathbb{I}}&=0,\quad \forall\,\boldsymbol{\omega}^h \in [V^h]^2,
	\end{aligned}
\end{equation}
while the update of the numerical solution remains given by~\eqref{eq3.8} and the curvature constraint by~\eqref{eqka}.

The following lemma reveals a key property of the reformulated scheme.

\begin{lemma}\label{lem4.1}
	Let the parameter interval $\mathbb{I}=[0,1]$ be uniformly partitioned with nodes $\rho_j=j/N$, $j=0,\dots,N$. Define the edge vectors $\mathbf{h}_j^{m,i}=\mathbf{X}^{m,i}(\rho_j)-\mathbf{X}^{m,i}(\rho_{j-1})$ and adopt the periodic convention $\mathbf h_{N+1}^{m,i}=\mathbf h_1^{m,i}$. For each $i=1,\dots,q$, suppose that $\bigl(\mathbf{X}^{m,i},\kappa^{m,i}\bigr)\in[V^h]^2\times V^h$ satisfies the second equation of~\eqref{eq4.4},
	\begin{equation}\label{eq4.5}
		\bigl( \kappa^{m,i},\; -(\partial_\rho\mathbf{X}^{m,i})^\perp \cdot \boldsymbol{\omega}^h \bigr)_{\mathbb{I}}^h - \frac{1}{L^{m,i}}\bigl( \partial_\rho\mathbf{X}^{m,i},\; \partial_\rho\boldsymbol{\omega}^h \bigr)_{\mathbb{I}}=0,\quad \forall\,\boldsymbol{\omega}^h \in [V^h]^2,
	\end{equation}
	with $L^{m,i}>0$. Then all edges of the stage curve have equal length, i.e.,
	\begin{equation}
		|\mathbf{h}_{j+1}^{m,i}|=|\mathbf{h}_{j}^{m,i}|,\quad j=1,\dots,N.
	\end{equation}
\end{lemma}

\begin{proof}
	Fix a node $\rho_j$ and choose the test function $\boldsymbol{\omega}^h=\varphi_j(\mathbf{h}_j^{m,i}+\mathbf{h}_{j+1}^{m,i})$, where $\varphi_j$ is the standard periodic nodal basis function. By the definition~\eqref{Newproduct} of the mass‑lumped inner product, the first term in~\eqref{eq4.5} becomes
	\begin{equation}\label{eq4.7}
		\bigl( \kappa^{m,i},\; -(\partial_\rho\mathbf{X}^{m,i})^\perp \cdot \boldsymbol{\omega}^h \bigr)_{\mathbb{I}}^h=-\frac{\kappa^{m,i}(\rho_j)}{2}\left(\mathbf{h}_j^{m,i}+\mathbf{h}_{j+1}^{m,i} \right)^\perp\cdot\left(\mathbf{h}_j^{m,i}+\mathbf{h}_{j+1}^{m,i}\right)=0,
	\end{equation}
	since $(\mathbf{v}+\mathbf{w})^\perp\cdot(\mathbf{v}+\mathbf{w})=0$ for any vectors $\mathbf{v},\mathbf{w}$.
	
	For the second term, uniform partitioning gives $\Delta\rho_j=1/N$ and $\partial_\rho\mathbf{X}^{m,i}=N\mathbf{h}_j^{m,i}$ on $I_j$. Hence
	\begin{equation}\label{eq4.8}
		\begin{aligned}
			\frac{1}{L^{m,i}}\bigl( \partial_\rho\mathbf{X}^{m,i},\; \partial_\rho\boldsymbol{\omega}^h \bigr)_{\mathbb{I}}&=\frac{N}{L^{m,i}}\left(\mathbf{h}_j^{m,i}-\mathbf{h}_{j+1}^{m,i} \right)\cdot\left(\mathbf{h}_j^{m,i}+\mathbf{h}_{j+1}^{m,i}\right)\\
			&=\frac{N}{L^{m,i}}\left(|\mathbf{h}_j^{m,i}|^2-|\mathbf{h}_{j+1}^{m,i}|^2 \right).
		\end{aligned}
	\end{equation}
	Substituting~\eqref{eq4.7} and~\eqref{eq4.8} into~\eqref{eq4.5} and using $L^{m,i}>0$ yields $|\mathbf{h}_{j+1}^{m,i}|=|\mathbf{h}_{j}^{m,i}|$, completing the proof.
\end{proof}

Based on Lemma~\ref{lem4.1}, we next show that solution of the reformulated coupled system \eqref{eq4.4} also satisfies the original stage formulation \eqref{eq3.7}.

\begin{theorem}\label{thm4.1}
	Let $\bigl(\mathbf{X}^{m,i},\kappa^{m,i}\bigr)\in[V^h]^2\times V^h$, $i=1,\dots,q$, be the stage solution of the coupled system~\eqref{eq4.4}. Then, for each $i$, the corresponding variables also satisfy the original formulation~\eqref{eq3.7}.
\end{theorem}
\begin{proof}
	For each stage $i$ and each element $I_j$, Lemma~\ref{lem4.1} gives $|\mathbf h_j^{m,i}|=L^{m,i}/N$. Since the partition is uniform, $\Delta\rho_j=1/N$, and therefore
	\begin{equation}
		|\partial_\rho\mathbf{X}^{m,i}|=\frac{|\mathbf{h}_{j}^{m,i}|}{\Delta\rho_j}=L^{m,i}, \quad \rho\in I_j.
	\end{equation}
	Thus the local metric factor in the original formulation is replaced by the global curve length $L^{m,i}$, which implies that the reformulated system~\eqref{eq4.4} coincides exactly with the original stage system~\eqref{eq3.7}.
\end{proof}

In analogy with the equivalence established in~\cite{2011BGN} between the reformulated and original first‑order schemes, we now prove the equivalence between~\eqref{eq4.4} and~\eqref{eq3.7} under a mild condition.

\begin{theorem}\label{thm4.2}
	Assume that $\mathbb{I}$ is uniformly partitioned and that the stage curves contain no degenerate elements. If the solution of the original RK stage scheme~\eqref{eq3.7} satisfies
	\begin{equation}\label{eq4.10}
		\frac{\mathbf{h}_{j+1}^{m,i}}{|\mathbf{h}_{j+1}^{m,i}|} \neq \frac{\mathbf{h}_j^{m,i}}{|\mathbf{h}_j^{m,i}|}, \quad j=1,\dots, N,\quad i=1,\dots,q,
	\end{equation}
	i.e., adjacent unit tangent vectors are distinct, then for each stage $i$ the scheme~\eqref{eq4.4} is equivalent to scheme~\eqref{eq3.7}.
\end{theorem}

\begin{proof}
	The implication from the reformulated system~\eqref{eq4.4} to the original system~\eqref{eq3.7} follows from Theorem~\ref{thm4.1}. To prove the converse, let $\mathbf{t}_j^{m,i}={\mathbf{h}_{j}^{m,i}}/{|\mathbf{h}_{j}^{m,i}|}$ for $j=1,\dots,N$, with $\mathbf{t}_{N+1}^{m,i}=\mathbf{t}_{1}^{m,i}$. Choosing the test function $\boldsymbol{\omega}^h=\varphi_j(\mathbf{h}_j^{m,i}+\mathbf{h}_{j+1}^{m,i})$ in the second equation of~\eqref{eq3.7} yields
	\begin{equation}
		\begin{aligned}
			\bigl( \partial_s^{m,i}\mathbf{X}^{m,i},\; \partial_s^{m,i}\boldsymbol{\omega}^h \bigr)_{\Gamma^{m,i}}&=\left(\mathbf{t}_j^{m,i}-\mathbf{t}_{j+1}^{m,i} \right)\cdot\left(\mathbf{h}_j^{m,i}+\mathbf{h}_{j+1}^{m,i}\right)\\
			&=\left(|\mathbf{h}_j^{m,i}|-|\mathbf{h}_{j+1}^{m,i}| \right)\left(1-\mathbf{t}_j^{m,i}\cdot\mathbf{t}_{j+1}^{m,i}\right).
		\end{aligned}
	\end{equation}
	As in~\eqref{eq4.7}, the first term of that equation vanishes. Consequently,
	\begin{equation}
		\left(|\mathbf{h}_j^{m,i}|-|\mathbf{h}_{j+1}^{m,i}| \right)\left(1-\mathbf{t}_j^{m,i}\cdot\mathbf{t}_{j+1}^{m,i}\right)=0.
	\end{equation}
	Condition~\eqref{eq4.10} ensures $1-\mathbf{t}_j^{m,i}\cdot\mathbf{t}_{j+1}^{m,i} \neq0$, which forces $|\mathbf{h}_{j}^{m,i}|=|\mathbf{h}_{j+1}^{m,i}|$. The uniform parameter partition then implies $|\partial_\rho\mathbf{X}^{m,i}|=L^{m,i}$, showing that any solution of~\eqref{eq3.7} also solves the reformulated scheme~\eqref{eq4.4}. This completes the proof.
\end{proof}

\begin{rek}
	The intermediate stages $\mathbf{X}^{m,i}$ produced by the reformulation~\eqref{eq4.4} satisfy mesh equidistribution (Lemma~\ref{lem4.1}). However, this does not imply that the numerical solution $\mathbf{X}^{m+1}$ at the next time level is equidistributed, because the RK update~\eqref{eq3.8} uses a weighted combination of stage velocities.
\end{rek}
\subsection{The damped Newton method}
For the coupled nonlinear system~\eqref{eq4.4}, we employ a monolithic damped Newton iteration that computes all RK stage solutions simultaneously. At the $k$‑th iteration, given current approximations $\mathbf X_k^{m,i}$ and $\kappa_k^{m,i}$ ($i=1,\dots,q$), we seek the Newton increments $\delta\mathbf X_k^{m,i}$ and $\delta\kappa_k^{m,i}$.

For convenience, define the stage velocities and curve lengths at the current iteration by
\begin{equation}\nonumber
\mathbf\Phi_k^{m,i}
=
\frac{1}{\tau}
\sum_{j=1}^{q}
\widehat a_{ij}
\left(\mathbf X_k^{m,j}-\mathbf X^m\right),
\qquad
L_k^{m,i}
=
\int_{\mathbb I}
\left|\partial_\rho\mathbf X_k^{m,i}\right|\,\mathrm{d}\rho.
\end{equation}
Their directional derivatives with respect to the increments are
\begin{equation}\nonumber
	\delta\mathbf\Phi_k^{m,i}
	=
	\frac{1}{\tau}\sum_{j=1}^{q}\widehat a_{ij}\delta\mathbf X_k^{m,j},
	\qquad
	\delta L_k^{m,i}
	=
	\int_{\mathbb I}
	\frac{\partial_\rho\mathbf X_k^{m,i}}
	{|\partial_\rho\mathbf X_k^{m,i}|}
	\cdot\partial_\rho\delta\mathbf X_k^{m,i}\,\mathrm{d}\rho.
\end{equation}

Linearizing~\eqref{eq4.4} with respect to all stage unknowns gives the Newton system: for each $i=1,\dots,q$, find $\delta\mathbf X_k^{m,i}\in[V^h]^2$ and $\delta\kappa_k^{m,i}\in V^h$ such that

\begin{subequations}\label{eq4.13}
	\begin{align}
		&
		\begin{aligned}
			&
			\left(
			\delta\mathbf\Phi_k^{m,i},
			-(\partial_\rho\mathbf X_k^{m,i})^\perp\varphi^h
			\right)_{\mathbb I}^h
			+
			\left(
			\mathbf\Phi_k^{m,i},
			-(\partial_\rho\delta\mathbf X_k^{m,i})^\perp
			\varphi^h
			\right)_{\mathbb I}^h
			\\
			&\quad+
			\frac{1}{L_k^{m,i}}
			\left(
			\partial_\rho\delta\kappa_k^{m,i},
			\partial_\rho\varphi^h
			\right)_{\mathbb I}
			-
			\frac{\delta L_k^{m,i}}{(L_k^{m,i})^2}
			\left(
			\partial_\rho\kappa_k^{m,i},
			\partial_\rho\varphi^h
			\right)_{\mathbb I}
			\\
			&=
			-
			\left(
			\mathbf\Phi_k^{m,i},
			-(\partial_\rho\mathbf X_k^{m,i})^\perp\varphi^h
			\right)_{\mathbb I}^h
			-
			\frac{1}{L_k^{m,i}}
			\left(
			\partial_\rho\kappa_k^{m,i},
			\partial_\rho\varphi^h
			\right)_{\mathbb I},
			\quad \forall\,\varphi^h\in V^h.
		\end{aligned}
		\label{eq4.13a}
		\\[3mm]
		&
		\begin{aligned}
			&
			\left(
			\delta\kappa_k^{m,i},
			-(\partial_\rho\mathbf X_k^{m,i})^\perp
			\cdot\boldsymbol\omega^h
			\right)_{\mathbb I}^h
			+
			\left(
			\kappa_k^{m,i},
			-(\partial_\rho\delta\mathbf X_k^{m,i})^\perp
			\cdot\boldsymbol\omega^h
			\right)_{\mathbb I}^h
			\\
			&\quad-
			\frac{1}{L_k^{m,i}}
			\left(
			\partial_\rho\delta\mathbf X_k^{m,i},
			\partial_\rho\boldsymbol\omega^h
			\right)_{\mathbb I}
			+
			\frac{\delta L_k^{m,i}}{(L_k^{m,i})^2}
			\left(
			\partial_\rho\mathbf X_k^{m,i},
			\partial_\rho\boldsymbol\omega^h
			\right)_{\mathbb I}
			\\
			&=
			-
			\left(
			\kappa_k^{m,i},
			-(\partial_\rho\mathbf X_k^{m,i})^\perp
			\cdot\boldsymbol\omega^h
			\right)_{\mathbb I}^h
			+
			\frac{1}{L_k^{m,i}}
			\left(
			\partial_\rho\mathbf X_k^{m,i},
			\partial_\rho\boldsymbol\omega^h
			\right)_{\mathbb I},
			\quad \forall\,\boldsymbol\omega^h\in[V^h]^2.
		\end{aligned}
		\label{eq4.13b}
	\end{align}
\end{subequations}
The equations are coupled through the stage velocity increments $\delta\mathbf\Phi_k^{m,i}$, which depend on the position increments at all $q$ stages. Hence the Newton direction is obtained by solving the entire linear system as a single block.

Once the increments are computed, all stage solutions are updated simultaneously:
\[
\begin{aligned}
	\mathbf X_{k+1}^{m,i}
	=
	\mathbf X_k^{m,i}
	+\alpha_k\delta\mathbf X_k^{m,i},\qquad
	\kappa_{k+1}^{m,i}
	=
	\kappa_k^{m,i}
	+\alpha_k\delta\kappa_k^{m,i},
	\qquad i=1,\ldots,q,
\end{aligned}
\]
where $\alpha_k\in(0,1]$ is a damping parameter determined by backtracking line search to ensure a sufficient decrease of the residual.

To monitor the convergence of the iteration, let  
$$\mathbf Z_k
=\left(
\mathbf X_k^{m,1},\kappa_k^{m,1},
\ldots,
\mathbf X_k^{m,q},\kappa_k^{m,q}
\right)^T
\in\mathbb R^{3Nq}$$ 
collect the finite‑element coefficient vectors of all stage unknowns, and let $\mathcal R(\mathbf Z_k)\in\mathbb R^{3Nq}$ denote the assembled nonlinear residual of~\eqref{eq4.13}. The iteration is terminated when
\[
\left\|\mathcal R(\mathbf Z_k)\right\|_2
\leq \mathrm{tol},
\]
with a prescribed tolerance $\mathrm{tol}=10^{-10}$. Denote by $k_*$ the iteration index at which this criterion is first met. The final stage solutions are then
\[
\mathbf X^{m,i}=\mathbf X_{k_*}^{m,i},\quad
\kappa^{m,i}=\kappa_{k_*}^{m,i},\quad
\mathbf\Phi^{m,i}
=\frac{1}{\tau}\sum_{j=1}^q\widehat a_{ij}
(\mathbf X_{k_*}^{m,j}-\mathbf X^m).
\]
The curve at the next time level is updated via
\[
\mathbf X^{m+1}
=\mathbf X^m+\tau\sum_{i=1}^q b_i\mathbf\Phi^{m,i}.
\]
In practice, a small diagonal regularization is added to the Newton Jacobian to enhance the robustness of the linear solver.

\begin{rek}
	Initial guesses for the stage positions and curvatures are constructed using a semi‑implicit BGN1 predictor~\cite{2007SISCBGN,2007JCPBGN} evaluated at the RK stage times $t_m+c_i\tau$. If needed, extrapolation from the previous time level or from the previously accepted stage velocity can also serve as alternative initializations.
\end{rek}

\begin{rek}
	Direct application of Newton's method to the original fully implicit RK system often suffers from convergence failures due to geometric nonlinearity. By performing Newton iterations on the equivalent reformulation of Barrett, Garcke, and N\"urnberg~\cite{2011BGN}—which replaces the local metric factor by the global curve length—we obtain a significantly more robust nonlinear solver while preserving exact equivalence with the original discrete scheme.
\end{rek}

\section{Numerical experiments}\label{Section5}

In this section, we present numerical experiments that validate the accuracy, area‑preservation, and computational efficiency of the proposed symplectic RK schemes for surface diffusion flow of curves. All simulations were conducted using FreeFem++ 4.14, and the results were post‑processed and visualized with MATLAB R2018a. The numerical methods compared in the following are summarized below:
\begin{itemize}
	\item AP‑GRK1: Area‑preserving symplectic RK scheme with coefficients \eqref{eqGRK1}.
	\item AP‑GRK2: Area‑preserving symplectic RK scheme with coefficients \eqref{eqGRK2}.
	\item AP‑GRK3: Area‑preserving symplectic RK scheme with coefficients \eqref{eqGRK3}.
	\item BGN1 \cite{2007SISCBGN,2007JCPBGN}: The classic first‑order Euler scheme.
	\item BGN‑BDF2 \cite{2024SISCBDFk}: The second‑order backward differentiation formula scheme.
	\item BGN‑PC \cite{2025JCPPC}: The second‑order predictor‑corrector scheme.
\end{itemize}

\begin{table}[tbhp]
	\caption{Manifold errors and observed convergence orders of the 1‑stage, second‑order AP‑GRK1 scheme at different time instants, with initial data \eqref{eq5.2} and the parameter relation $\tau=0.05h$.}
	\label{GRK1C}
	\centering
	\renewcommand{\arraystretch}{1.4}
	\begin{tabular}{|c|cc|cc|cc|}
		\hline
		\multicolumn{7}{|c|}{AP‑GRK1 }\\ \hline
		\multirow{2}{*}{$\tau$}
		& \multicolumn{2}{c|}{$T=0.1$}
		& \multicolumn{2}{c|}{$T=0.5$}
		& \multicolumn{2}{c|}{$T=1$} \\
		\cline{2-7}
		& $\mathcal{E}_{\mathcal{M}}$ & Order
		& $\mathcal{E}_{\mathcal{M}}$ & Order
		& $\mathcal{E}_{\mathcal{M}}$ & Order \\ \hline
		$1/400$
		& 1.0507E-1 & --
		& 1.1009E-1 & --
		& 9.1758E-2 & -- \\ \hline
		$1/800$
		& 2.6442E-2 & 1.9904
		& 2.8502E-2 & 1.9495
		& 2.3466E-2 & 1.9673 \\ \hline
		$1/1600$
		& 6.6281E-3 & 1.9962
		& 7.2191E-3 & 1.9812
		& 5.9246E-3 & 1.9858 \\ \hline
		$1/3200$
		& 1.6567E-3 & 2.0003
		& 1.8142E-3 & 1.9925
		& 1.5180E-3 & 1.9645 \\ \hline
	\end{tabular}
\end{table}

\subsection{Accuracy verification}

When tangential mesh motion is present, errors measured with respect to a specific parametrization are not meaningful. Because the BGN formulation permits tangential redistribution of mesh nodes, we adopt the reparametrization‑invariant manifold distance introduced in \cite{2021IMAZhao} and also used in \cite{2024JCPCN}. For two closed polygonal curves $\Gamma_1$ and $\Gamma_2$, this distance is defined by
\begin{equation}
	\mathcal{M}(\Gamma_1,\Gamma_2)
	:= |(\Omega_1\setminus\Omega_2)\cup(\Omega_2\setminus\Omega_1)|
	=|\Omega_1|+|\Omega_2|-2|\Omega_1\cap\Omega_2|,
\end{equation}
where $\Omega_1$ and $\Omega_2$ denote the bounded regions enclosed by $\Gamma_1$ and $\Gamma_2$, respectively.

\begin{table}[tbhp]
	\caption{Manifold errors and observed convergence orders of the 2‑stage, fourth‑order AP‑GRK2 scheme at different time instants, with initial data \eqref{eq5.2} and the parameter relation $\tau=0.05h^{1/2}$.}
	\label{GRK2C}
	\centering
	\renewcommand{\arraystretch}{1.4}
	\begin{tabular}{|c|cc|cc|cc|}
		\hline
		\multicolumn{7}{|c|}{AP‑GRK2 }\\ \hline
		\multirow{2}{*}{$\tau$}
		& \multicolumn{2}{c|}{$T=0.1$}
		& \multicolumn{2}{c|}{$T=0.5$}
		& \multicolumn{2}{c|}{$T=1$} \\
		\cline{2-7}
		& $\mathcal{E}_{\mathcal{M}}$ & Order
		& $\mathcal{E}_{\mathcal{M}}$ & Order
		& $\mathcal{E}_{\mathcal{M}}$ & Order \\ \hline
		$1/80$
		& 2.0311E-1 & --
		& 2.0775E-1 & --
		& 1.7733E-1 & -- \\ \hline
		$1/160$
		& 1.2834E-2 & 3.9842
		& 1.3820E-2 & 3.9100
		& 1.1502E-2 & 3.9465 \\ \hline
		$1/320$
		& 8.0246E-4 & 3.9994
		& 8.6763E-4 & 3.9935
		& 7.1868E-4 & 4.0004 \\ \hline
		$1/640$
		& 5.0192E-5 & 3.9989
		& 5.4236E-5 & 3.9998
		& 4.5351E-5 & 3.9861 \\ \hline
	\end{tabular}
\end{table}

In the absence of an exact solution, the convergence order is estimated by comparing numerical solutions obtained on successive space‑time refinement levels \cite{21SINUMBao,2025SISCLagrange,2026JiangarXiv}. The error $\mathcal{E}_{\mathcal{M}}$ and the corresponding convergence order are defined as
\begin{equation}\nonumber
	\mathcal{E}_{\mathcal{M}}(T,\tau_1,\tau_2):=\mathcal{M}\left(\Gamma_{h_1}^{T/\tau_1},\Gamma_{h_2}^{T/\tau_2}\right),
	\quad
	\mathrm{Order}
	:=\log\left(\frac{\mathcal{E}_{\mathcal{M}}(T,\tau_1,\tau_2)}
	{\mathcal{E}_{\mathcal{M}}(T,\tau_2,\tau_3)}\right)
	\bigg/
	\log\left(\frac{\tau_1}{\tau_2}\right).
\end{equation}
Here, $\Gamma_h^{T/\tau}:=\mathbf{X}_h^{T/\tau}(\mathbb{I})$, and the time step $\tau_j$ and mesh size $h_j$ are related by $\tau_j=0.05(h_j)^{2/p}$, where $p$ denotes the formal order of the numerical method. In this work we consider $p=2,4,6$.

We test the convergence rates using the initially parametrized ellipse given in \cite{2025JCPPC}:
\begin{equation}\label{eq5.2}
	\mathbf{X}(\rho)=(2\cos (2\pi\rho),\sin (2\pi\rho)),\quad \rho\in\mathbb{I}=[0,1].
\end{equation}

Tables \ref{GRK1C}--\ref{GRK3C} indicate that the AP‑GRK1, AP‑GRK2, and AP‑GRK3 schemes achieve the expected second‑, fourth‑, and sixth‑order convergence rates, respectively, at all three time instants $T=0.1,\,0.5,\,1$. This confirms the well‑known fact that a $q$‑stage Gauss RK method is of order $2q$ \cite{2010FengKang,Hairer}. This provides a new route for developing high-order time integrators within the BGN framework.

\begin{table}[tbhp]
	\caption{Manifold errors and observed convergence orders of the 3‑stage, sixth‑order AP‑GRK3 scheme at different time instants, with initial data \eqref{eq5.2} and the parameter relation $\tau=0.05h^{1/3}$.}
	\label{GRK3C}
	\centering
	\renewcommand{\arraystretch}{1.4}
	\begin{tabular}{|c|cc|cc|c|cc|}
		\hline
		\multicolumn{8}{|c|}{AP‑GRK3}\\ \hline
		\multirow{2}{*}{$\tau$}
		& \multicolumn{2}{c|}{$T=0.1$}
		& \multicolumn{2}{c|}{$T=0.5$}
		& \multirow{2}{*}{$\tau$}& \multicolumn{2}{c|}{$T=1$} \\
		\cline{2-5}\cline{7-8}
		& $\mathcal{E}_{\mathcal{M}}$ & Order
		& $\mathcal{E}_{\mathcal{M}}$ & Order
		& & $\mathcal{E}_{\mathcal{M}}$ & Order \\ \hline
		$1/40$
		& 7.556E-1 & --
		& 7.377E-1 & --
		& $1/40$ & 6.243E-1 & -- \\ \hline
		$1/80$
		& 1.345E-2 & 5.812
		& 1.448E-2 & 5.671
		& $1/60$ & 6.220E-2  & 5.688 \\ \hline
		$1/160$
		& 2.102E-4 & 6.000
		& 2.272E-4 & 5.994
		& $1/90$ & 5.550E-3 & 5.960 \\ \hline
		$1/320$
		& 3.364E-6 & 5.965
		& 3.524E-6 & 6.011
		& $1/135$ & 4.858E-4 & 6.007 \\ \hline
	\end{tabular}
\end{table}

\begin{figure}[!http]
	\centering
	\begin{minipage}[!ht]{0.5\linewidth}
		\centering
		\subfloat[Normalized perimeter vs.\ step]{\includegraphics[width=6.4cm,height=4.5cm]{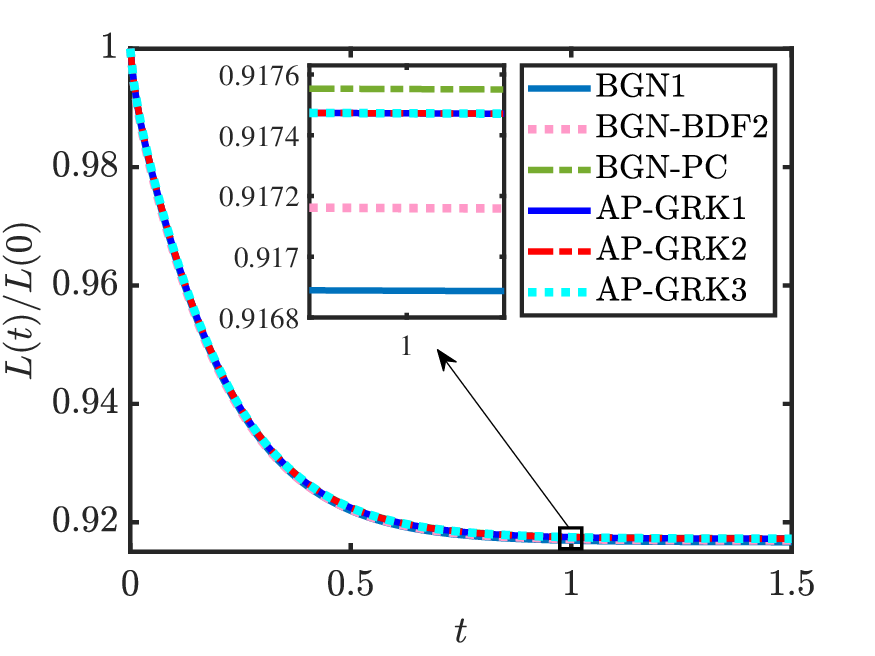}\label{Coma}}
	\end{minipage}%
	\begin{minipage}[!ht]{0.5\linewidth}
		\centering
		\subfloat[Relative area loss vs.\ step]{\includegraphics[width=6.4cm,height=4.5cm]{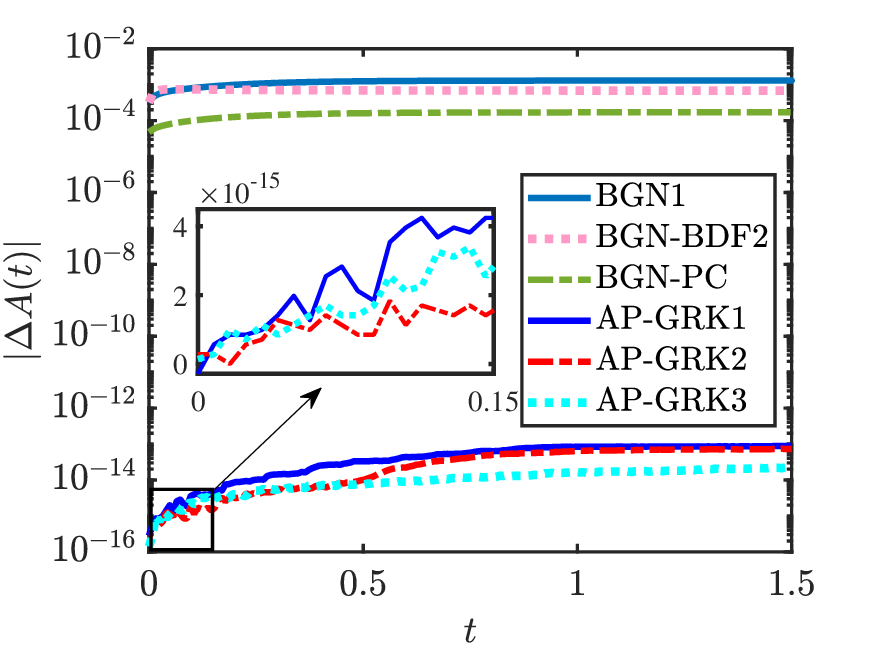}\label{Comb}}
	\end{minipage}\\
	\begin{minipage}[!ht]{0.5\linewidth}
		\centering
		\subfloat[Mesh ratio vs.\ step]{\includegraphics[width=6.4cm,height=4.5cm]{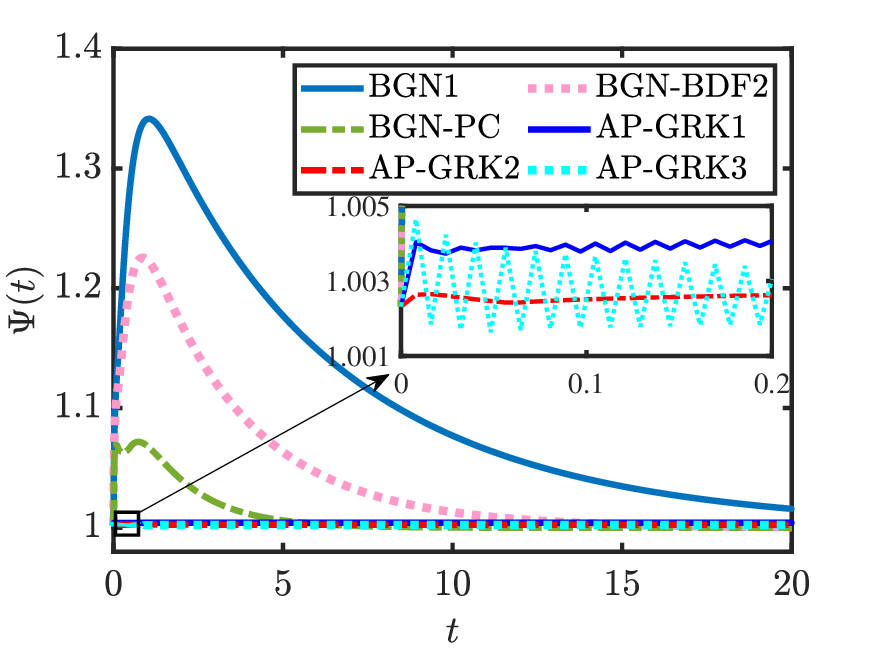}\label{Comc}}
	\end{minipage}%
	\begin{minipage}[!ht]{0.5\linewidth}
		\centering
		\subfloat[Iteration count per step]{\includegraphics[width=6.4cm,height=4.5cm]{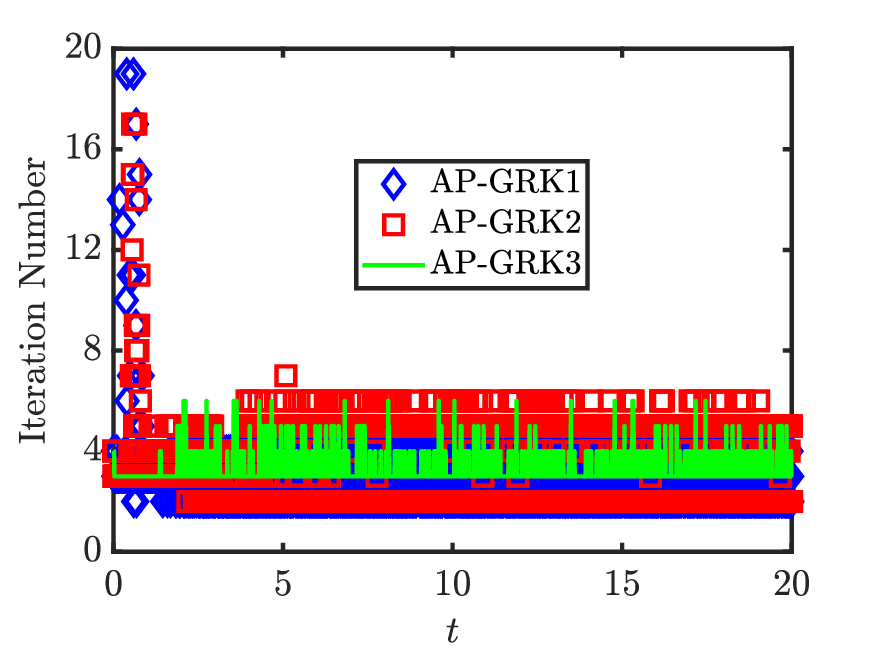}\label{Comd}}
	\end{minipage}
	\vspace{-3mm}
	\caption{Evolution of geometric quantities and iteration counts for initial data \eqref{eq5.2}, with $\tau=1/124$ and $N=80$.}
	\label{ComPicture}   
\end{figure}

\subsection{Comparison of advantages}

To highlight the benefits of the proposed symplectic RK schemes, we compare them against three classical BGN‑type methods: BGN1 \cite{2007JCPBGN}, BGN‑BDF2 \cite{2024SISCBDFk}, and BGN‑PC \cite{2025JCPPC}.

We assess stability by monitoring three diagnostic quantities during the evolution: the relative area loss $\Delta A(t)$, the normalized perimeter $L(t)/L(0)$, and the mesh‑distribution measure $\Psi(t)$. These are defined as
\begin{equation}\nonumber
	\left.\Delta A(t)\right|_{t=t_m}
	=\frac{A^m-A^0}{A^0},
	\quad
	\left.\frac{L(t)}{L(0)}\right|_{t=t_m}
	=\frac{L^m}{L^0},
	\quad
	\left.\Psi(t)\right|_{t=t_m}
	=\frac{\max_j |\mathbf{h}_j^m|}{\min_j |\mathbf{h}_j^m|},\quad m\geq0,
\end{equation}
where $A^m$ and $L^m$ denote the area enclosed by and the perimeter of the polygonal curve $\mathbf{X}^m$, respectively.

Using initial data \eqref{eq5.2} with $\tau=1/124$ and $N=80$, Figure~\ref{ComPicture} compares the geometric evolution produced by the different schemes. The following observations are made:
\begin{enumerate}[label=\arabic*)]
	\item All AP‑GRK schemes and the other BGN variants exhibit a monotonically decreasing perimeter (Figure~\ref{Coma}).

	\item The AP‑GRK schemes preserve area exactly (Figure~\ref{Comb}), whereas BGN1, BGN‑BDF2, and BGN‑PC incur area loss. This numerically confirms Theorem~\ref{thmAP}.

	\item For the ellipse evolution, the AP‑GRK schemes maintain a mesh ratio $\Psi(t)\rightarrow 1$ from the outset (Figure~\ref{Comc}), indicating good mesh quality. In contrast, the other  BGN schemes require a long transient before equidistribution is achieved.

	\item Although the AP‑GRK schemes are fully implicit and couple multiple intermediate stages, the number of Newton iterations drops rapidly to fewer than 8 per step (Figure~\ref{Comd}), demonstrating the efficiency of the tailored solver from Section~\ref{Section4}.
\end{enumerate}

\begin{figure}[!http]
	\centering
	\begin{minipage}[!ht]{0.33\linewidth}
		\centering
		\subfloat[Initial curve]{\includegraphics[width=4cm,height=3cm]{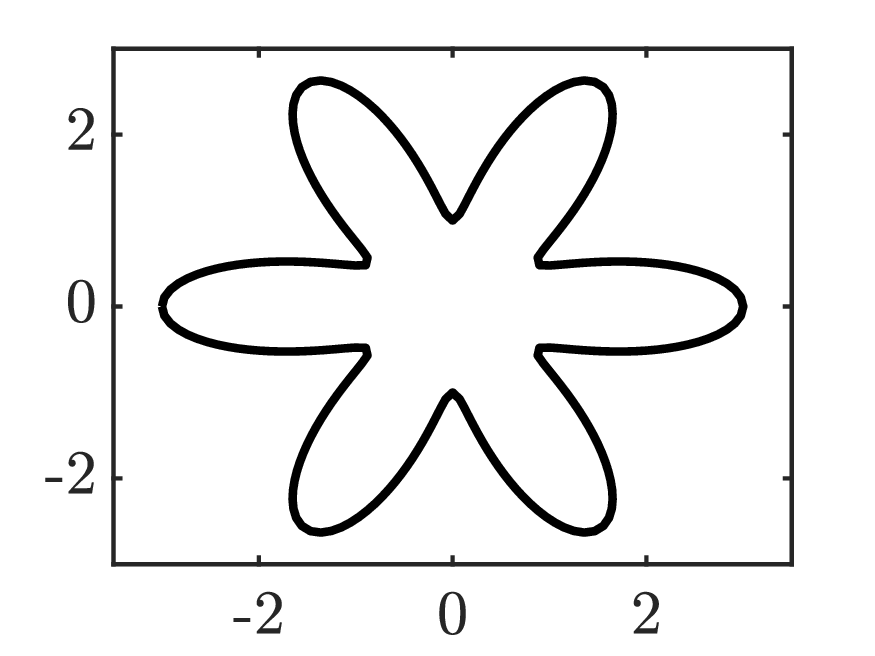}}
	\end{minipage}%
	\begin{minipage}[!ht]{0.33\linewidth}
		\centering
		\subfloat[$t=0.00111$]{\includegraphics[width=4cm,height=3cm]{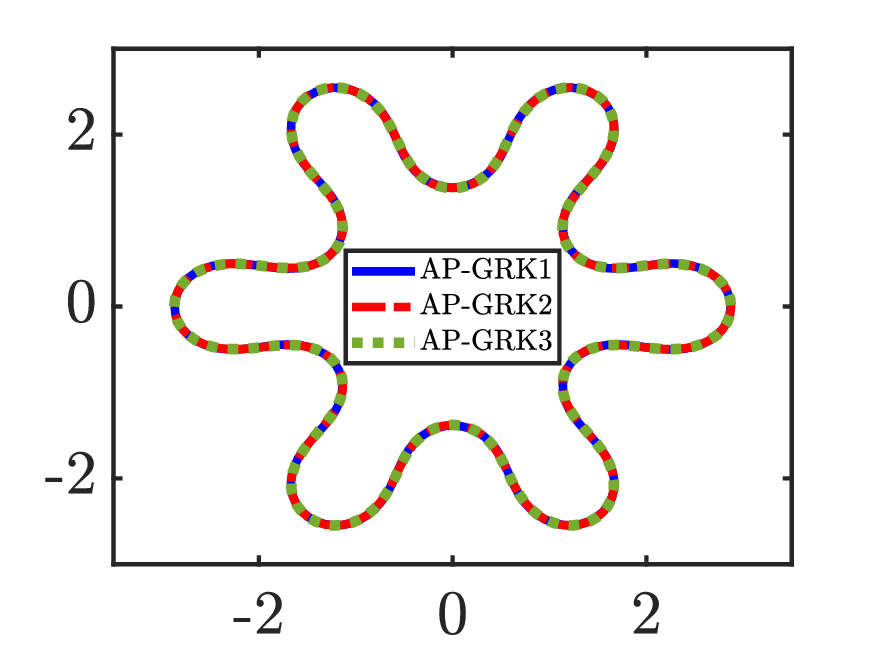}}
	\end{minipage}%
	\begin{minipage}[!ht]{0.33\linewidth}
		\centering
		\subfloat[$t=0.0333$]{\includegraphics[width=4cm,height=3cm]{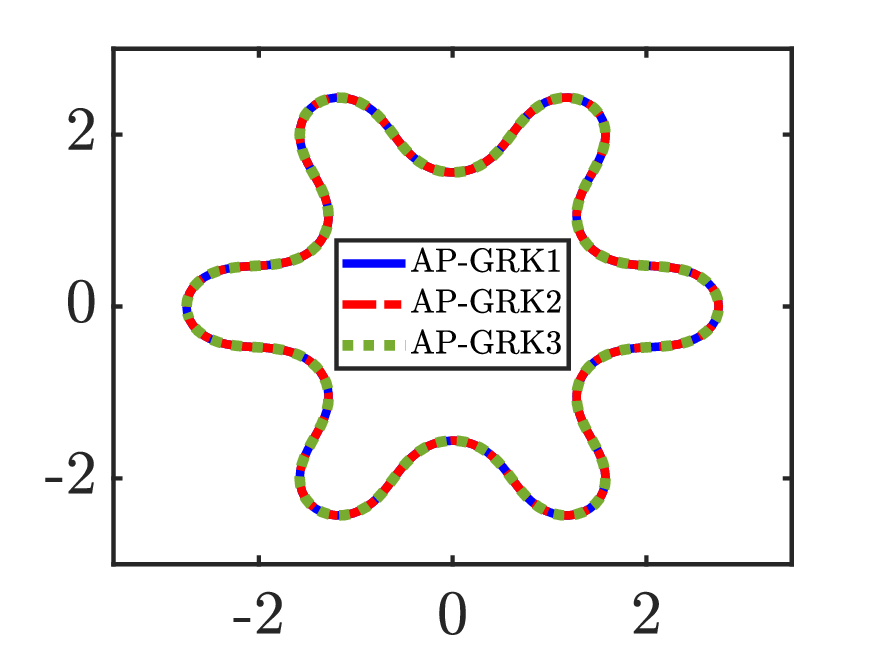}}
	\end{minipage}\\
	\begin{minipage}[!ht]{0.33\linewidth}
		\centering
		\subfloat[$t=0.0667$]{\includegraphics[width=4cm,height=3cm]{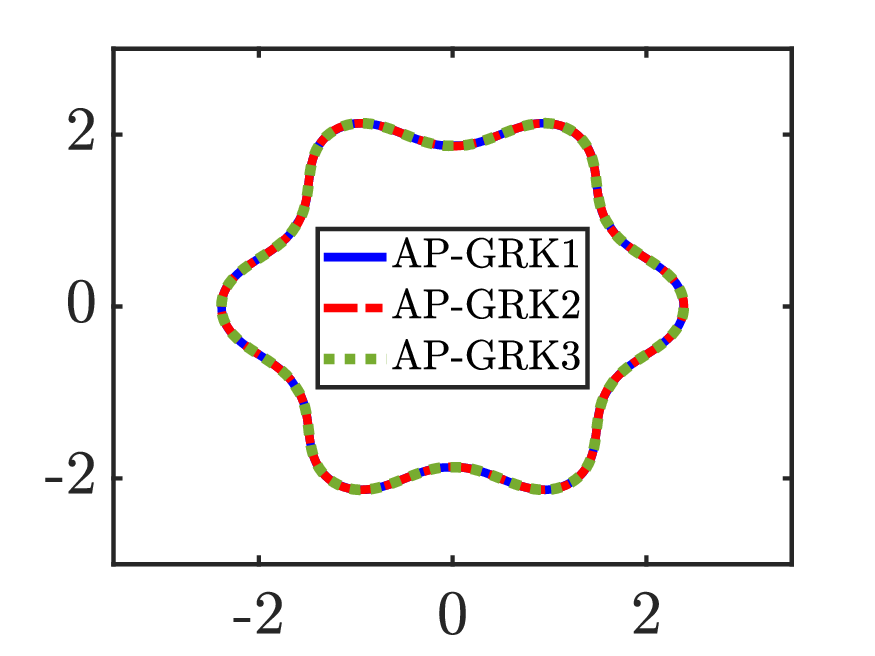}}
	\end{minipage}%
	\begin{minipage}[!ht]{0.33\linewidth}
		\centering
		\subfloat[$t=0.0889$]{\includegraphics[width=4cm,height=3cm]{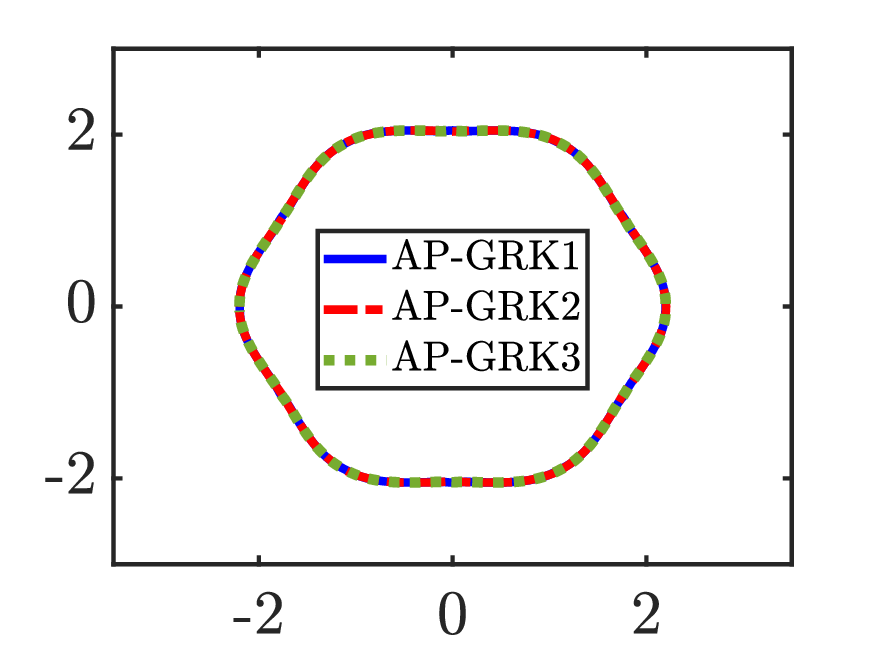}}
	\end{minipage}%
	\begin{minipage}[!ht]{0.33\linewidth}
		\centering
		\subfloat[$t=0.5$]{\includegraphics[width=4cm,height=3cm]{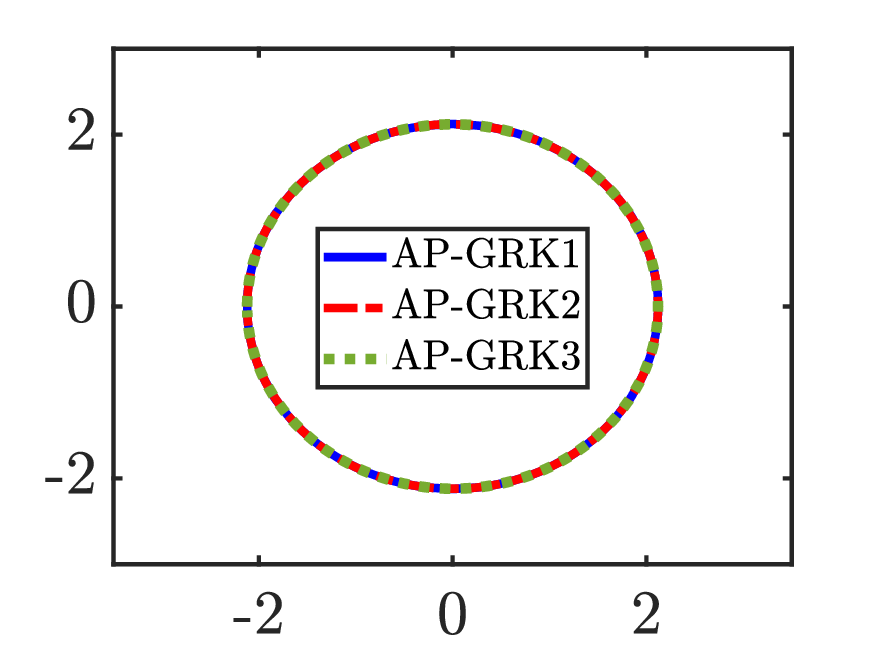}}
	\end{minipage}
	\vspace{-3mm}
	\caption{Evolution of an initially flower‑shaped curve \eqref{Flower} at different times with $\tau=1/180$ and $N=256$.}
	\label{EFlower}   
\end{figure}

\begin{figure}[!http]
	\centering
	\begin{minipage}[!ht]{0.5\linewidth}
		\centering
		\subfloat[Normalized perimeter]{\includegraphics[width=6.4cm,height=4.5cm]{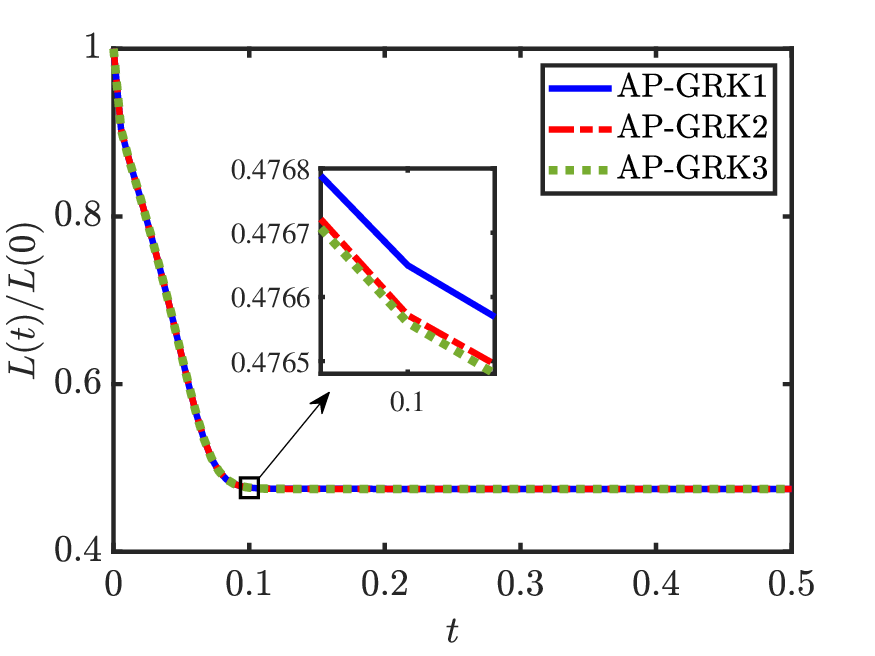}}
	\end{minipage}%
	\begin{minipage}[!ht]{0.5\linewidth}
		\centering
		\subfloat[Relative area loss]{\includegraphics[width=6.4cm,height=4.5cm]{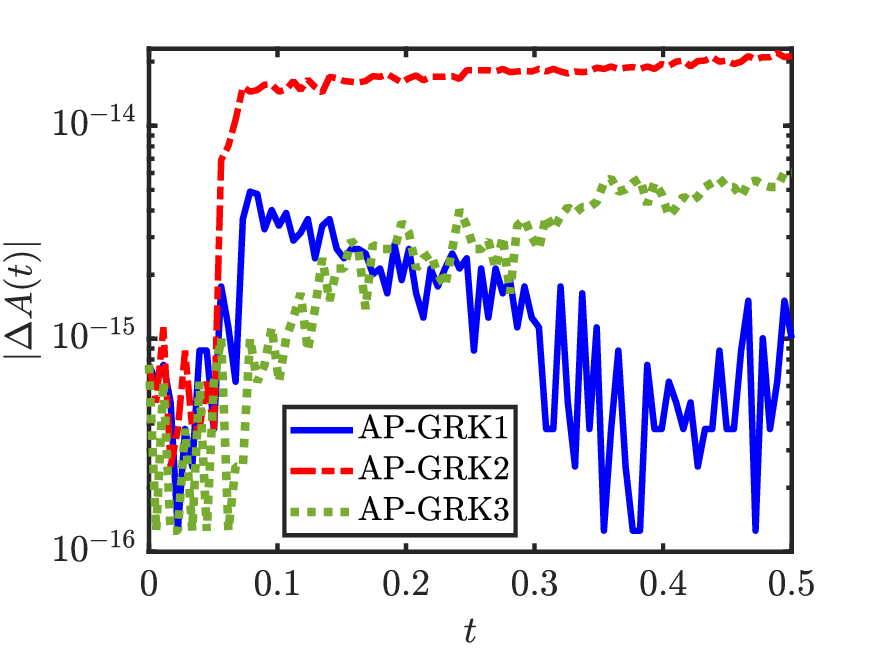}}
	\end{minipage}\\
	\begin{minipage}[!ht]{0.5\linewidth}
		\centering
		\subfloat[Mesh ratio]{\includegraphics[width=6.4cm,height=4.5cm]{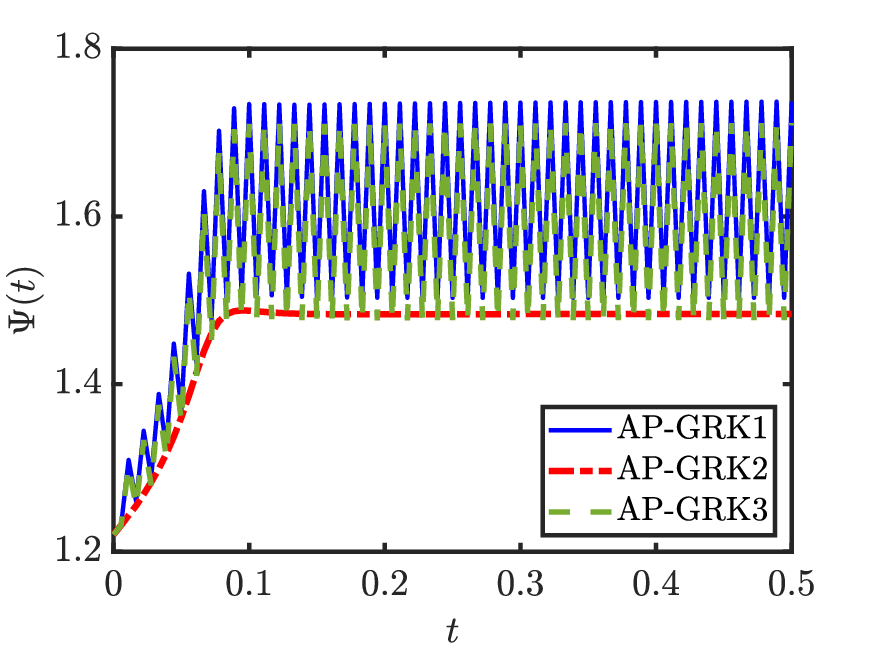}}
	\end{minipage}%
	\begin{minipage}[!ht]{0.5\linewidth}
		\centering
		\subfloat[Iteration count]{\includegraphics[width=6.4cm,height=4.5cm]{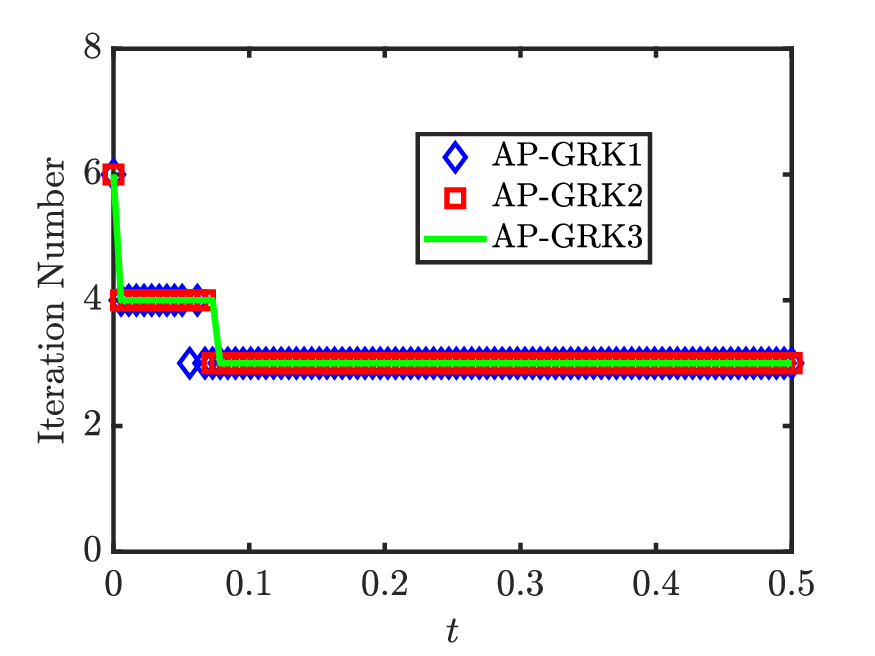}}
	\end{minipage}
	\vspace{-3mm}
	\caption{Evolution of geometric quantities and iteration counts for initial data \eqref{Flower}, with $\tau=1/180$ and $N=256$.}
	\label{GEFlower}   
\end{figure}

\subsection{Evolution of closed curves}

In this subsection, we demonstrate the capability of the AP‑GRK schemes to capture the smooth evolution of complex curves toward circular equilibria. For challenging initial shapes, we first perform a single step with the structure‑preserving Bao‑Zhao method \cite{21SINUMBao} to obtain well‑behaved initial data. This pre‑processing yields a favorable tangent‑vector distribution without introducing area loss, thereby enhancing robustness. The following initial curves are considered:
\begin{enumerate}[label=\roman*)]
	\item Flower‑shaped curve with six petals \cite{21SINUMBao,2025JCPPC}:
	\begin{equation}\label{Flower}
		\begin{cases}
			x=[2+\cos(12\pi\rho)]\cos(2\pi\rho),\\
			y=[2+\cos(12\pi\rho)]\sin(2\pi\rho),
		\end{cases}
		\quad \rho\in\mathbb{I}=[0,1].
	\end{equation}
	\item Non‑convex Mikula‑\v{S}ev\v{c}ovi\v{c} curve \cite{2001MS}:
	\begin{equation}\label{Ms}
		\begin{cases}
			x=\cos(2\pi\rho),\qquad \rho\in\mathbb{I}=[0,1],\\
			y=\frac{1}{2}\left( \sin(\cos(2\pi\rho))+\sin(2\pi\rho)[0.7+\sin(2\pi\rho)\sin^2(6\pi\rho)]\right).
		\end{cases}
	\end{equation}
	\item A cross‑shaped curve formed by two rectangles of aspect ratio 4:1 intersecting perpendicularly at their centers.
\end{enumerate}

\begin{figure}[!http]
	\centering
	\begin{minipage}[!ht]{0.33\linewidth}
		\centering
		\subfloat[Initial curve]{\includegraphics[width=4cm,height=3cm]{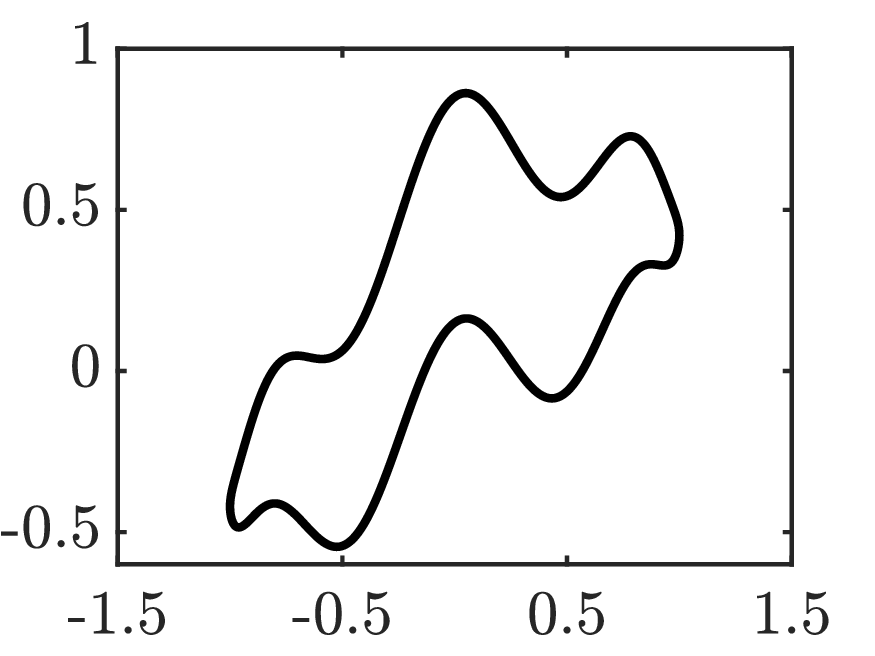}}
	\end{minipage}%
	\begin{minipage}[!ht]{0.33\linewidth}
		\centering
		\subfloat[$t=0.00256$]{\includegraphics[width=4cm,height=3cm]{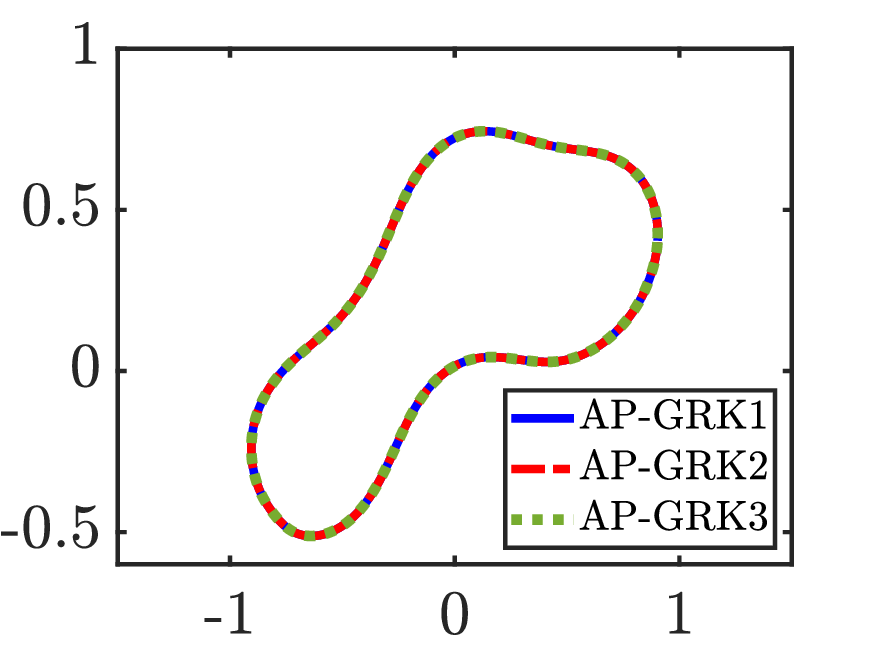}}
	\end{minipage}%
	\begin{minipage}[!ht]{0.33\linewidth}
		\centering
		\subfloat[$t=0.00513$]{\includegraphics[width=4cm,height=3cm]{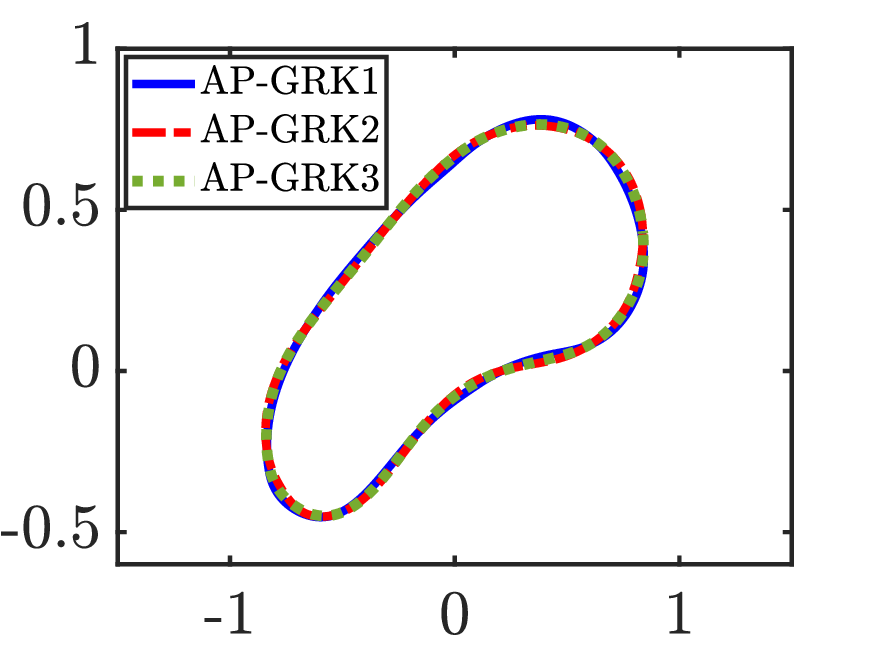}}
	\end{minipage}\\
	\begin{minipage}[!ht]{0.33\linewidth}
		\centering
		\subfloat[$t=0.0103$]{\includegraphics[width=4cm,height=3cm]{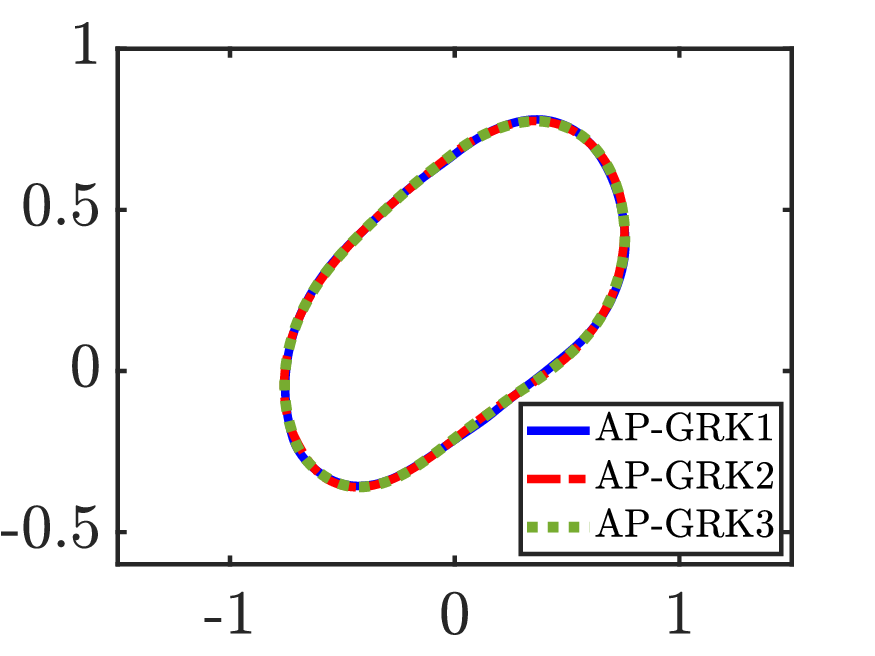}}
	\end{minipage}%
	\begin{minipage}[!ht]{0.33\linewidth}
		\centering
		\subfloat[$t=0.0205$]{\includegraphics[width=4cm,height=3cm]{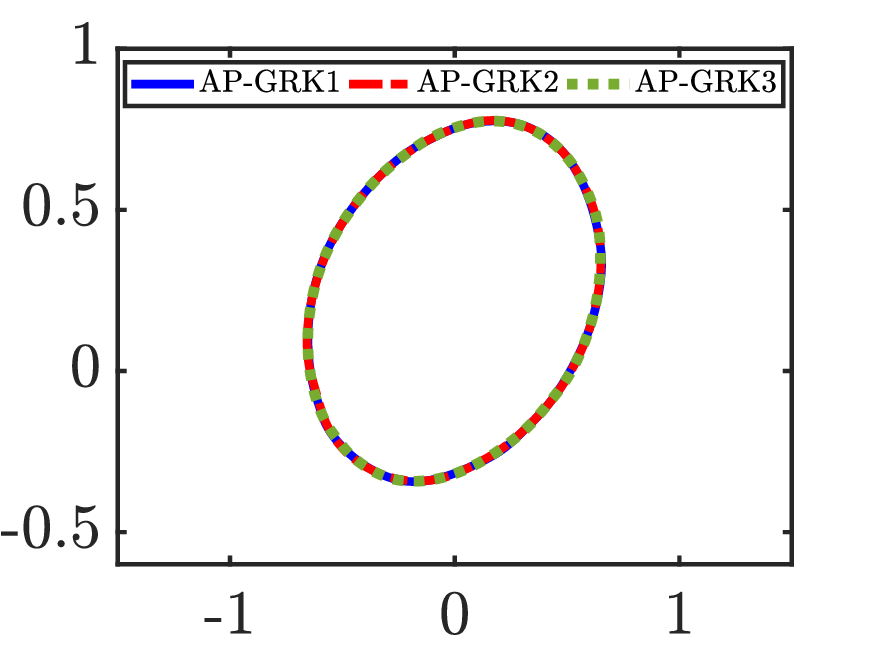}}
	\end{minipage}%
	\begin{minipage}[!ht]{0.33\linewidth}
		\centering
		\subfloat[$t=0.1$]{\includegraphics[width=4cm,height=3cm]{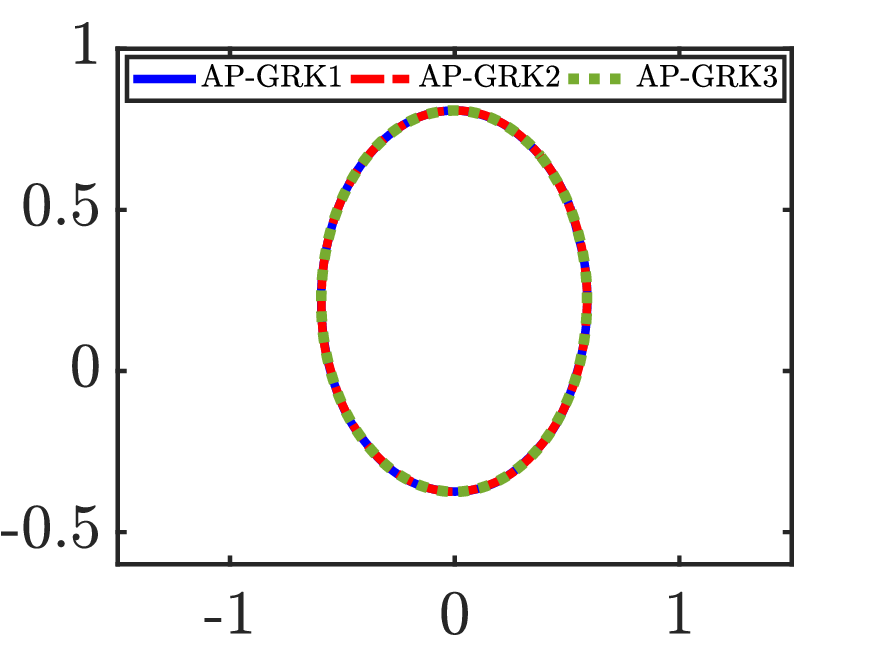}}
	\end{minipage}
	\vspace{-3mm}
	\caption{Evolution of an initially Mikula‑\v{S}ev\v{c}ovi\v{c} curve \eqref{Ms} at different times with $\tau=1/390$ and $N=1024$.}
	\label{EMs}   
\end{figure}

\begin{figure}[!http]
	\centering
	\begin{minipage}[!ht]{0.5\linewidth}
		\centering
		\subfloat[Normalized perimeter]{\includegraphics[width=6.4cm,height=4.5cm]{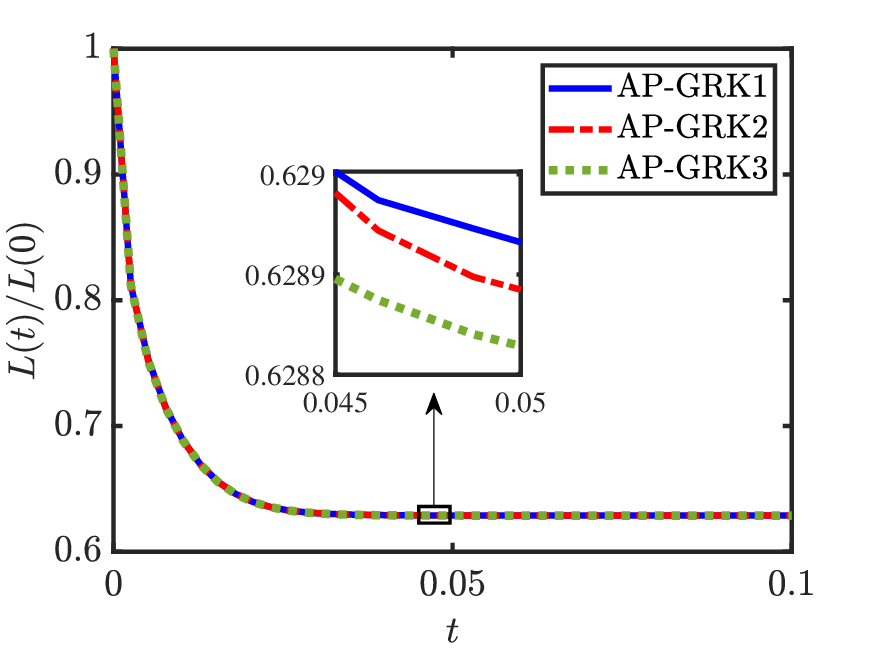}}
	\end{minipage}%
	\begin{minipage}[!ht]{0.5\linewidth}
		\centering
		\subfloat[Relative area loss]{\includegraphics[width=6.4cm,height=4.5cm]{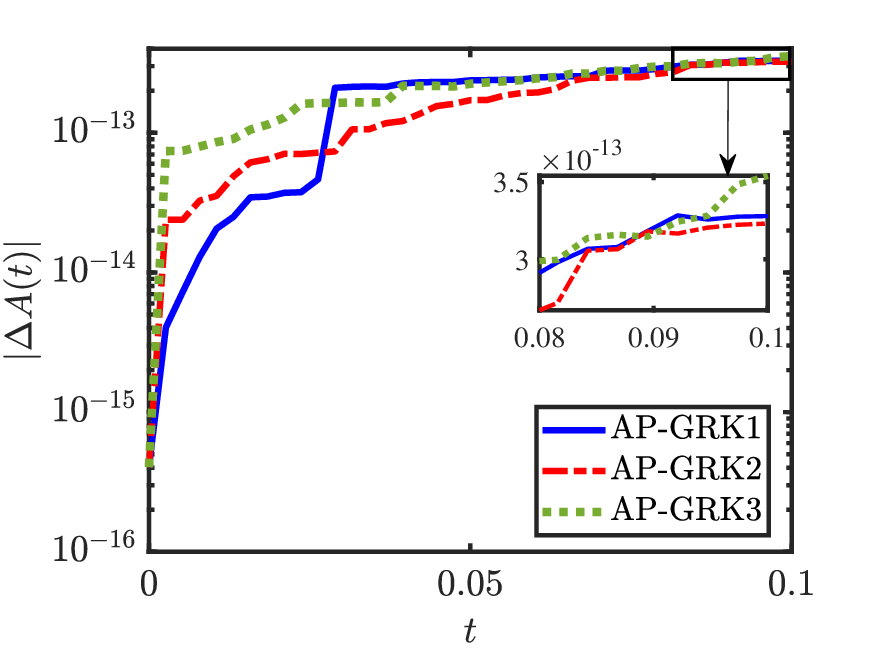}}
	\end{minipage}\\
	\begin{minipage}[!ht]{0.5\linewidth}
		\centering
		\subfloat[Mesh ratio]{\includegraphics[width=6.4cm,height=4.5cm]{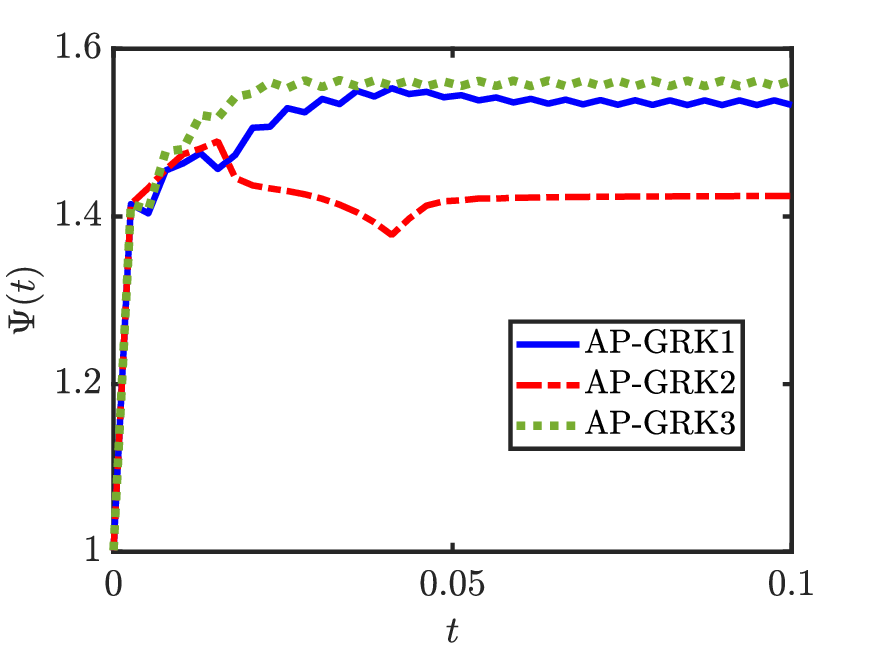}}
	\end{minipage}%
	\begin{minipage}[!ht]{0.5\linewidth}
		\centering
		\subfloat[Iteration count]{\includegraphics[width=6.4cm,height=4.5cm]{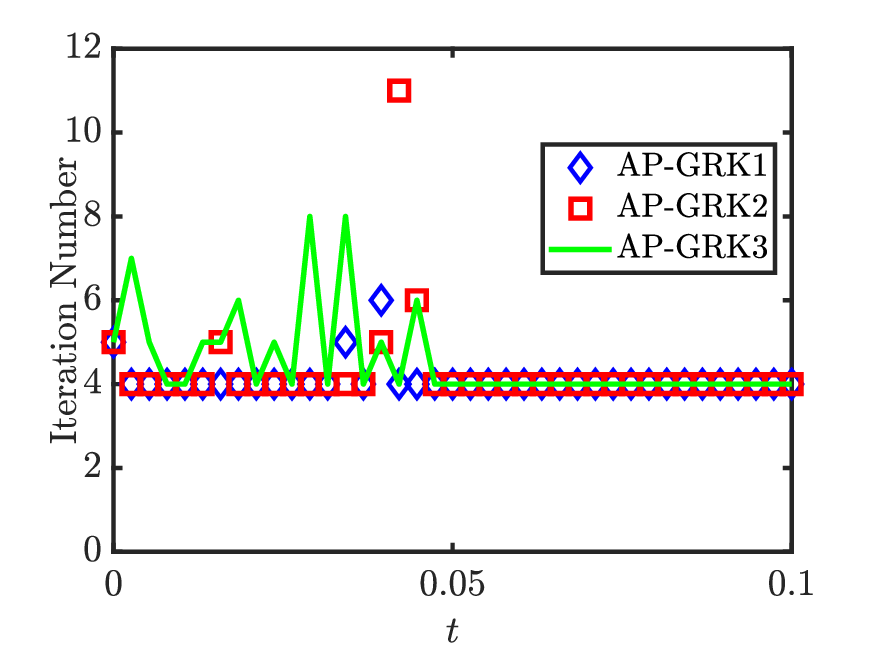}\label{MsIN}}
	\end{minipage}
	\vspace{-3mm}
	\caption{Evolution of geometric quantities and iteration counts for initial data \eqref{Ms}, with $\tau=1/390$ and $N=1024$.}
	\label{GEMs}   
\end{figure}

For the flower‑shaped curve \eqref{Flower} with $\tau=1/180$ and $N=256$, Figure~\ref{EFlower} shows that all three AP‑GRK schemes accurately capture the evolution: the six petals gradually fade and the curve eventually approaches a circle. Figure~\ref{GEFlower} confirms that area is conserved exactly, perimeter decays monotonically, mesh quality remains high, and only a small number of iterations are required.

For the Mikula‑\v{S}ev\v{c}ovi\v{c} curve \eqref{Ms} with $\tau=1/390$ and $N=1024$, Figure~\ref{EMs} illustrates the evolution computed by the AP‑GRK schemes: sharp corners are progressively smoothed and the curve tends toward a circular shape. The corresponding geometric diagnostics in Figure~\ref{GEMs} again verify area conservation, perimeter decay, and good mesh quality, while the iteration counts remain modest.

\begin{figure}[!http]
	\centering
	\begin{minipage}[!ht]{0.5\linewidth}
		\centering
		\subfloat[Evolution of cross‑shaped curve]{\includegraphics[width=6.4cm,height=4.5cm]{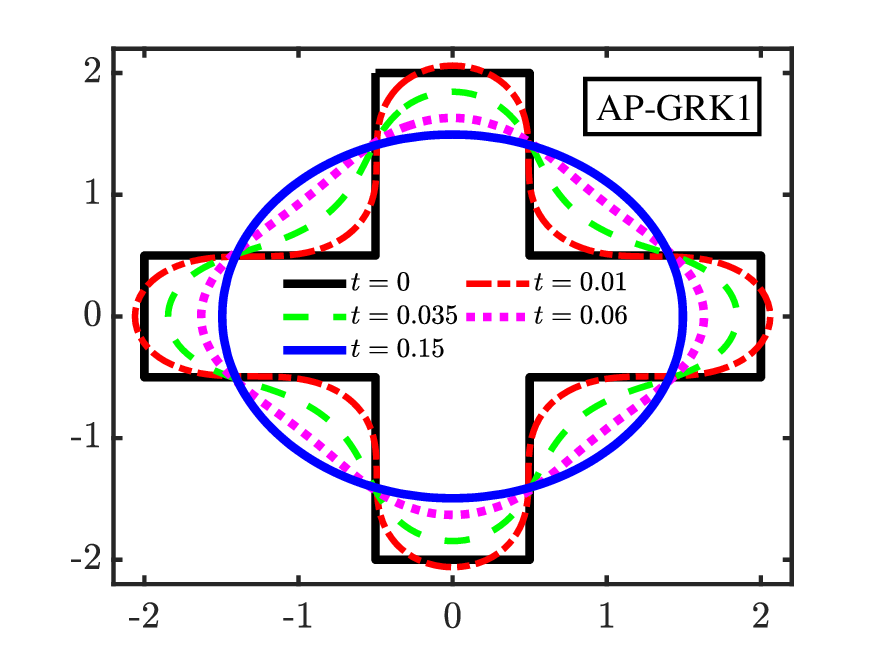}}
	\end{minipage}%
	\begin{minipage}[!ht]{0.5\linewidth}
		\centering
		\subfloat[Normalized perimeter and iterations]{\includegraphics[width=6.4cm,height=4.5cm]{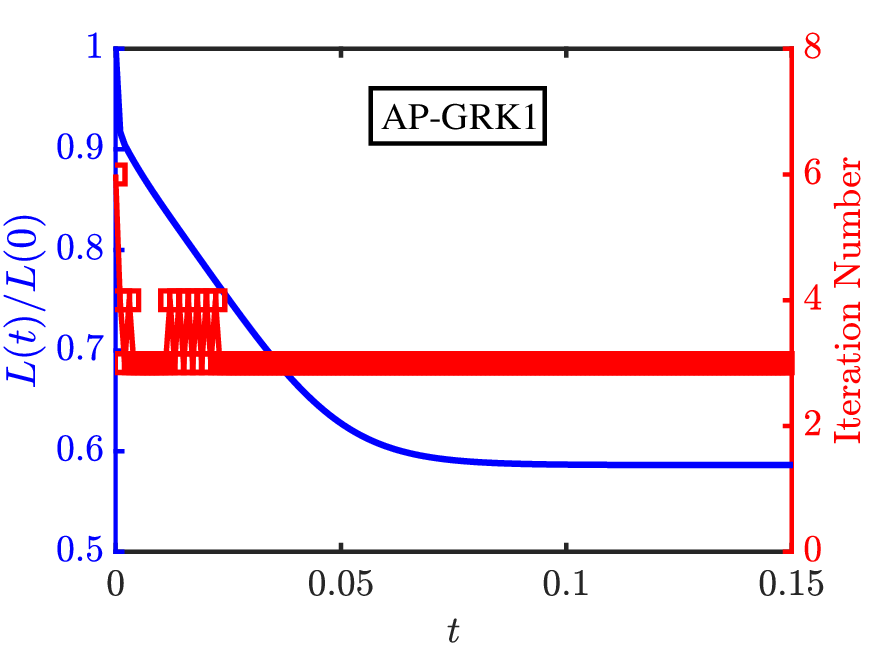}}
	\end{minipage}\\
	\begin{minipage}[!ht]{0.5\linewidth}
		\centering
		\subfloat[Mesh ratio]{\includegraphics[width=6.4cm,height=4.5cm]{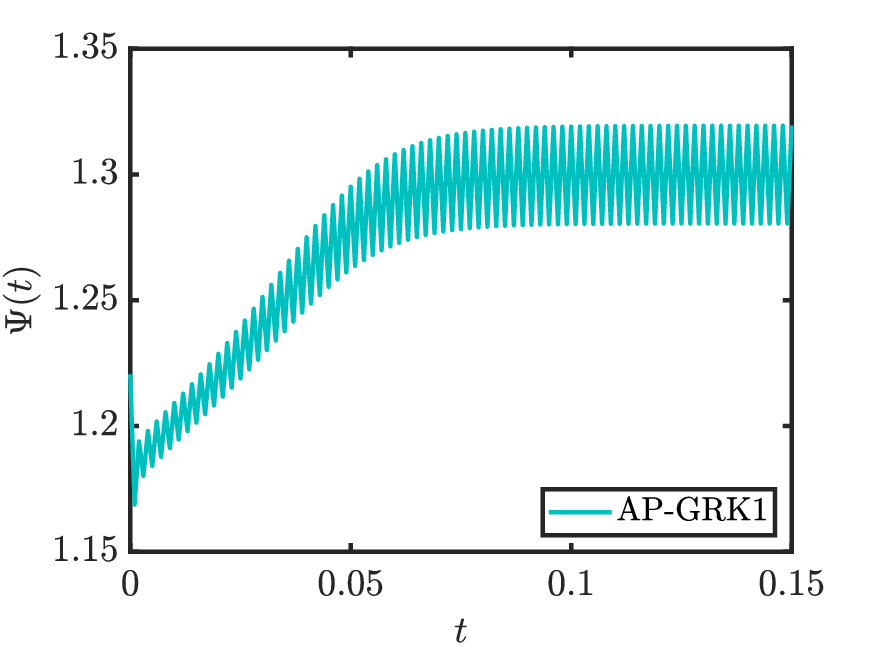}}
	\end{minipage}%
	\begin{minipage}[!ht]{0.5\linewidth}
		\centering
		\subfloat[Relative area loss]{\includegraphics[width=6.4cm,height=4.5cm]{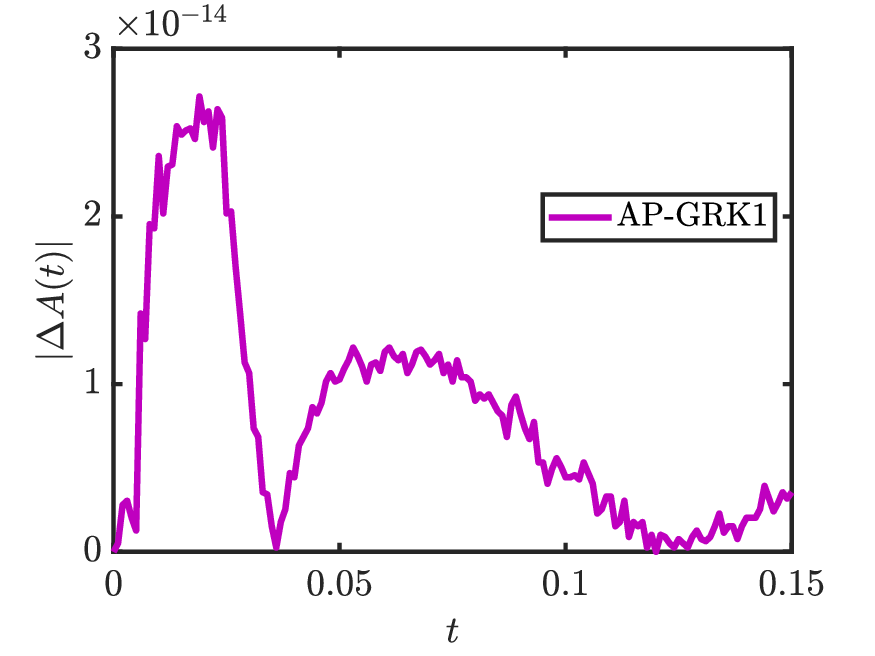}}
	\end{minipage}
	\vspace{-3mm}
	\caption{Evolution of the cross‑shaped curve and its geometric quantities using the AP‑GRK1 scheme, with $\tau=1/1000$ and $N=256$.}
	\label{ECsc}   
\end{figure}

Finally, we simulate the evolution of the cross‑shaped curve using the AP‑GRK1 scheme with $\tau=1/1000$ and $N=256$. The results are presented in Figure~\ref{ECsc}. Right angles gradually become rounded and concave regions expand outward until a circular shape is attained. As in previous examples, AP‑GRK1 preserves area exactly, exhibits perimeter decay, and maintains excellent mesh quality, further validating the accuracy and efficiency of the proposed methods.

\section{Conclusions}\label{Section6}
We developed high‑order area‑preserving RK  schemes for surface diffusion flow of curves by combining the parametric finite element method (BGN formulation) with symplectic RK time discretisation. This approach simultaneously delivers unconditional area preservation, high‑order temporal accuracy, and efficient solution of the fully implicit system. We proved  that symplectic RK schemes unconditionally preserve enclosed area due to their intrinsic symplectic structure, without introducing additional constraints,   auxiliary variables,  or Lagrange multipliers to ensure structure preservation. Adopting $q$-stage Gauss RK methods (symplectic, order $2q$) yields a family of schemes with tunable high‑order accuracy (second, fourth, and sixth) while strictly maintaining exact area preservation. The resulting fully implicit and highly nonlinear system is addressed efficiently by means of an equivalent discrete formulation coupled with a monolithic damped Newton solver, numerical results  show the number of iterations per step rapidly drops to fewer than 8, demonstrating the robustness and efficiency of the tailored solver. Extensive tests validate area conservation, high‑order convergence rates, perimeter decay, and favourable mesh quality. The proposed approach provides a powerful tool for simulating surface diffusion flow, with future extensions to three‑dimensional surfaces, adaptive time‑stepping, and related curvature‑driven geometric flows.

\appendix
\appendixnotitle \label{appendixA}
The Butcher tableaux of the one‑ to three‑stage Gauss RK methods \cite{2010FengKang,Hairer} are listed below:
\begin{itemize}
	\item The 1‑stage, second‑order Gauss RK method:
	\begin{equation}\label{eqGRK1}
		\begin{array}{c|c}
			\mathbf c & \mathbf A\\ \hline
			&  \mathbf b^T
		\end{array}
		=
		\begin{array}{c|c}
			\frac12 & \frac12\\ \hline
			& 1
		\end{array}.
	\end{equation}
	\item The 2‑stage, fourth‑order Gauss RK method:
	\begin{equation}\label{eqGRK2}
		\begin{array}{c|c}
			\mathbf c & \mathbf A\\ \hline
			& \mathbf b^T
		\end{array}
		=
		\begin{array}{c|cc}
			\frac{3-\sqrt3}{6} & \frac14 & \frac14-\frac{\sqrt3}{6}\\[2mm]
			\frac{3+\sqrt3}{6} & \frac14+\frac{\sqrt3}{6} & \frac14\\ \hline
			& \frac12 & \frac12
		\end{array}.
	\end{equation}
	\item The 3‑stage, sixth‑order Gauss RK method:
	\begin{equation}\label{eqGRK3}
		\begin{array}{c|c}
			\mathbf c & \mathbf A\\ \hline
			& \mathbf b^T
		\end{array}
		=
		\begin{array}{c|ccc}
			\frac{5-\sqrt{15}}{10}
			& \frac{5}{36}
			& \frac{2}{9}-\frac{\sqrt{15}}{15}
			& \frac{5}{36}-\frac{\sqrt{15}}{30}\\[2mm]
			\frac12
			& \frac{5}{36}+\frac{\sqrt{15}}{24}
			& \frac{2}{9}
			& \frac{5}{36}-\frac{\sqrt{15}}{24}\\[2mm]
			\frac{5+\sqrt{15}}{10}
			& \frac{5}{36}+\frac{\sqrt{15}}{30}
			& \frac{2}{9}+\frac{\sqrt{15}}{15}
			& \frac{5}{36}\\ \hline
			& \frac{5}{18}
			& \frac{4}{9}
			& \frac{5}{18}
		\end{array}.
	\end{equation}
\end{itemize}
All of the above schemes satisfy the symplectic condition \eqref{Symplectic}, as can be readily verified.

\bibliographystyle{siamplain}
\bibliography{references}
\end{document}